\documentclass[12pt]{amsart}
\usepackage{enumerate,a4}
\usepackage{amssymb}

\allowdisplaybreaks

\theoremstyle{plain}
  \newtheorem{theorem}{Theorem}[section]
  \newtheorem{proposition}[theorem]{Proposition}
  \newtheorem{lemma}[theorem]{Lemma}
  \newtheorem{corollary}[theorem]{Corollary}
  \newtheorem{conjecture}[theorem]{Conjecture}
\theoremstyle{definition}

\numberwithin{equation}{section}

\newcommand{\integers}{\mathbb{Z}}

\newcommand{\link}{\mathrm{lk}}

\newcommand{\fvec}{\mathfrak{f}}
\newcommand{\hvec}{\mathfrak{h}}

\newcommand{\depth}{\mathrm{depth}}
\newcommand{\prdim}{\mathrm{pdim}}
\newcommand{\reg}{\mathrm{reg}}
\newcommand{\height}{\mathrm{height}}
\newcommand{\sd}{\mathrm{sd}}
\newcommand{\xx}{\mathbf{\underline{x}}}
\newcommand{\End}{\mathrm{end}}
\newcommand{\des}{\mathrm{des}}
\newcommand{\Hilb}{\mathrm{Hilb}}
\newcommand{\IM}{\mathrm{Im}}
\newcommand{\KER}{\mathrm{Ker}}

\begin{document}
  
\title{The Multiplicity Conjecture for Barycentric Subdivisions}
  
\author{Martina Kubitzke}
\address{Fachbereich Mathematik und Informatik\\
      Philipps-Universit\"at Marburg\\
      35032 Marburg, Germany}
\email{kubitzke@mathematik.uni-marburg.de}

\author{Volkmar Welker}
\address{Fachbereich Mathematik und Informatik\\
    Philipps-Universit\"at Marburg\\
    35032 Marburg, Germany}
\email{welker@mathematik.uni-marburg.de}
  
\thanks{}
  
\keywords{Barycentric subdivision, Stanley-Reisner ideal, minimal free resolution,
multiplicity}
  
\subjclass{ }
  
\begin{abstract}
  For a simplicial complex $\Delta$ we study the effect of barycentric
  subdivision on ring theoretic invariants of its Stanley-Reisner ring.
  In particular, for Stanley-Reisner rings of barycentric subdivisions we 
  verify a conjecture by Huneke and Herzog \& Srinivasan,
  that relates the multiplicity of a standard graded $k$-algebra to the 
  product of the maximal and minimal 
  shifts in its minimal free resolution up to the 
  height.
  On the way to proving the conjecture we develop new and list 
  well known results on behavior of dimension, Hilbert series, 
  multiplicity, local cohomology, depth and regularity when passing from
  the Stanley-Reisner ring of $\Delta$ to the one of its barycentric
  subdivision.
\end{abstract} 

\maketitle
  
\section{Introduction}
  For a simplicial complex $\Delta$ on ground set $\Omega$ its Stanley-Reisner 
  ideal $I_\Delta$ is the ideal in $S = k[x_\omega~|~\omega \in \Omega]$
  generated by the monomials $\xx_A := \prod_{\omega \in A} x_\omega$
  for $A \subseteq \Omega$ and $A \not\in \Delta$. Many combinatorial invariants
  of $\Delta$ are encoded in ring-theoretic invariants of its Stanley-Reisner 
  ring $k[\Delta] := S/I_\Delta$. Here we are interested in the behavior of
  these invariants when passing from $k[\Delta]$ to $k[\sd(\Delta)]$,
  where $\sd(\Delta)$ denotes the barycentric subdivision of $\Delta$.
  Recall, that $\sd(\Delta)$ is the simplicial complex on ground set 
  $\mathring{\Delta} := \Delta \setminus \{ \emptyset\}$ whose simplices
  are flags $A_0 \subsetneq A_1 \subsetneq \cdots \subsetneq A_i$ of 
  elements $A_j \in \mathring{\Delta}$, $0 \leq j \leq i$.
  Note, that throughout the paper we assume that if $\Delta$ is a simplicial
  complex on ground set $\Omega$ then $\{ \omega \} \in \Delta$ for all
  $\omega \in \Omega$. In particular, $I_\Delta$ will not contain any 
  variable. 

  Our main result is the verification of the multiplicity conjecture by
  Huneke and Herzog \& Srinivasan \cite{HS} for $k[\sd(\Delta)]$.
  In recent years, this conjecture has attracted attention
  from commutative algebra and combinatorics (see for example \cite{Fr}, \cite{G}, 
  \cite{GSS}, \cite{GT}, 
  \cite{HZ}, \cite{MNR1}, \cite{MNR2}, \cite{MR}, \cite{NS}, \cite{R}).
  Here we provide another link to combinatorics and add a large
  class of rings for which the conjecture holds.

  In general, for a standard graded $k$-algebra $A = T/I$, where $T = 
  k[x_1, \ldots, x_n]$, the conjecture relates the multiplicity of $A$ and 
  the shifts in the minimal free resolution of $A$ up to its height.  
  More precisely, let
  $$0 \rightarrow \bigoplus_{j \geq 0} T(-j)^{\beta_{r,j}} \rightarrow
  \cdots \rightarrow \bigoplus_{j \geq 0} T(-j)^{\beta_{1,j}}
  \rightarrow T \rightarrow A \rightarrow 0$$

  be the minimal free resolution of $A$ as a $T$-module. 
  Let $e(A)$ denote the multiplicity of $A$, $h = \height(I)$ be the height or 
  codimension of $I$ and set $M_i = \max \{~j~|~j \geq 0$ and $\beta_{i,j} \neq 0\}$ respectively
  $m_i=\min\{~j~|~j \geq 0$ and $\beta_{i,j} \neq 0\}$.
  Then the conjecture states:

  \begin{conjecture}[Multiplicity Conjecture] \label{MultiplicityConjecture}
     $$e(A) \leq \frac{1}{h!} \prod_{i = 1}^h M_i.$$
	 If $A$ is Cohen-Macaulay then 
	 $$e(A)\geq \frac{1}{h!}\prod_{i=1}^h m_i.$$
	 If $A$ is Cohen-Macaulay then equality holds if and only if
	 $A$ has a pure resolution over $T$.
  \end{conjecture}

  %There is also a conjectured lower bound in the Cohen-Macaulay case, which we
  %were not able to attack by our methods.
  We refer the reader for more background in commutative algebra to the
  books by Eisenbud \cite{E} and Bruns \& Herzog \cite{BH}.

  Thus our main result states:

  \begin{theorem} \label{MainResult}
    Let $\Delta$ be a simplicial complex. Then the
    Multiplicity Conjecture holds for $k[\sd(\Delta)]$.
  \end{theorem}

  In independent work Novik and Swartz \cite{NS} have verified the
  the multiplicity conjecture for some important classes of Cohen-Macaulay 
  simplicial complexes. The classes of simplicial complexes treated by
  Novik and Swartz and our Theorem \ref{MainResult} overlap in 
  a small fraction of either class. The lower bound in Theorem 
  \ref{MainResult} is also a consequence of parallel work by
  Michael Goff.

  For the proof of Theorem \ref{MainResult} we need to study the behavior
  of a few ring theoretic invariants when passing from $k[\Delta]$ to
  $k[\sd(\Delta)]$. We take this as an opportunity to list in Section
  \ref{RingInvariants} the relation of the most important ring theoretic
  invariants of $k[\Delta]$ and $k[\sd(\Delta)]$.

  The proof of Theorem \ref{MainResult} will then be given in Section
  \ref{ProofOfMainResult} and will rely on the Hochster formula for
  the Betti numbers of the minimal free resolution of a Stanley-Reisner
  ring $k[\Delta]$. 

\section{Invariants for Barycentric Subdivisions}\label{RingInvariants}
  \subsection{Basic Definitions}
    Before we can discuss the behavior of ring theoretic invariants when passing from $k[\Delta]$ 
    to $k[\sd(\Delta)]$, we need to introduce some basic notation about simplicial complexes.

    For the formulation of the results and proofs, we adopt the following standard notation for 
    simplicial complexes (see \cite{B} for more details). For $F\in \Delta$
    we denote by $\partial F$ the simplicial complex of all $G \subsetneq F$ 
    that lie in the boundary of the simplex $F$. We call an element
    $F$ of $\Delta$ a face of $\Delta$. An inclusionwise maximal face is called
    facet. For a face $F$ its dimension is given as $\dim F = |F|-1$ and the dimension of
    $\Delta$ is the maximum dimension of one of its faces. 
    The vector $\fvec^\Delta := (f_{-1}^\Delta, \ldots, f_{\dim \Delta}^\Delta)$ where 
    $f_i^\Delta$ is the number of $i$-dimensional faces of $\Delta$ is called the $f$-vector
    of $\Delta$. The vector $\hvec^\Delta := (h_0^\Delta, \ldots, h_{\dim \Delta +1}^\Delta)$
    defined by $$\sum_{i=0}^{\dim \Delta +1} h_i^\Delta t^{\dim \Delta -i +1} =
    \sum_{i=0}^{\dim \Delta+1} f_{i-1}^\Delta (t-1)^{\dim \Delta +1-i}$$ is called
    the $h$-vector of $\Delta$. 
    For a face $F  \in \Delta$ we write $\link_\Delta(F) := \{ G \in \Delta~|~
    F \cap G = \emptyset,~F \cup G \in \Delta \}$ for the link of $F$ in $\Delta$.  
    By $\widetilde{H}_i(\Delta;k)$ we denote the $i$-th reduced simplicial
    homology group with coefficients in $k$. Also we use $[n]$ to denote for a natural 
    number $n$ the set $\{ 1,\ldots, n\}$. 

   \subsection{Krull Dimension} \label{Krull-dimension}
     Since $\dim \Delta=\dim \sd(\Delta)$ it follows that 
     $$\dim k[\Delta]=\dim \Delta +1 = \dim \sd(\Delta)+1 = \dim k[\sd(\Delta)].$$

   \subsection{Hilbert Series} \label{Hilbert-Series}
     \begin{proposition}\cite[Theorem 2.2]{BW}
       Let $\Delta$ be a $(d-1)$-dimensional simplicial complex.
       Then:
       \begin{eqnarray*}
         \Hilb\;(k[\Delta],t)&=&\frac{h^{\Delta}_{0}+h^{\Delta}_{1}t+\cdots+h^{\Delta}_{d}t^d}{(1-t)^d}\\
         \Hilb\;(k[\sd(\Delta)],t)&=& \frac{h^{\sd(\Delta)}_{0}+h^{\sd(\Delta)}_{1}t+
                \cdots+h^{\sd(\Delta)}_{d}t^d}{(1-t)^d}\\
         &=&\frac{\displaystyle{\sum_{j=0}^{d}}
         \Big(\displaystyle{\sum_{i=0}^{d}}h_i^\Delta A(d+1,j,i+1)\Big)t^j}{(1-t)^d}
        \end{eqnarray*}
        where $A(d+1,j,i+1)$ denotes the number of permutations $\sigma \in S_{d+1}$ such that 
        $\sigma(1)=i+1$ and $\des(\sigma):=
              \#\left\{l \in [d]\hspace{5pt}|\hspace{5pt}\sigma(l)>\sigma(l+1)\right\}=j$.
      \end{proposition}

    \subsection{Local Cohomology} \label{Local-Cohomology}
      We denote by $H^i(k[\Delta])=H^i_{\mathfrak{m}}(k[\Delta])$ the $i$-th local cohomology 
      module of $k[\Delta]$ with respect to $\mathfrak{m}=\left(x_1,\ldots,x_n\right)$ where 
      $n=f_0^{\Delta}$ (for more background see \cite{BS}).
      We recall the $\integers$-graded version of a theorem of Hochster for the Hilbert series of 
      the $i$-th local cohomology of $k[\Delta]$.

      \begin{proposition}[see Theorem 5.3.8 in \cite{BH}]\label{Hochster} 
        Let $\Delta$ be a simplicial complex. Then the $\mathbb{Z}$-graded Hilbert series of the 
        $i$-th local cohomology module of $k[\Delta]$ is given by
        $$\Hilb\;\left(H^i\left (k[\Delta] \right),t \right)=\sum_{F \in \Delta} 
        \dim_k \widetilde{H}_{i-|F|-1} \left (\link_{\Delta} F;k\right )\cdot 
        \left ( \frac{1}{t-1} \right )^{|F|}.$$
      \end{proposition}

      We also need the following simple lemma about links in barycentric subdivisions, whose 
      verification is left to the reader.

      \begin{lemma}\label{lemma:no1}
        Let $\Delta$ be a simplicial complex, $\sd(\Delta)$ its barycentric subdivision. 
        Then for a face $F$ of $\Delta$ and a flag $F_1\subsetneq \ldots F_r:=F$ of 
        $\sd(\Delta)$ it holds that
        $$\widetilde{H}_{i-|F|-1}\left(\link_{\Delta}F;k\right)=\widetilde{H}_{i-r-1}
        \left(\link_{\sd(\Delta)}\left(F_1\subsetneq \ldots \subsetneq F_r\right);k\right).$$ 
      \end{lemma}

      \begin{proposition}
        Let $\Delta$ be a $(d-1)$-dimensional simplicial complex. Then the $\mathbb{Z}$-graded 
        Hilbert series of the $i$-th local cohomology module of $k[\sd(\Delta)]$ is given by
        $$\Hilb\;(H^i(k[\sd(\Delta)]),t) = \dim_k \widetilde{H}_{i-1} \left (\Delta; k \right) +
        $$

        $$\sum_{m=1}^{d} \sum_{\substack{F \in \Delta \\ |F|=m}}
          \left (\frac{\displaystyle{\sum_{k=0}^{m-1}} 
          \Big| \left\{ \sigma \in S_m \; |\; \des(\sigma)=k \right\} \Big| 
          \cdot t^k}{(t-1)^m}\right )
          \cdot \dim_k \widetilde{H}_{i-m-1}\left (\link_{\Delta} F;k \right ).
        $$
      \end{proposition}
      \begin{proof}
         For $F\in \mathring{\Delta}$ we set 
         $$\sd(\Delta)[F]:=\left\{F_1\subsetneq\ldots\subsetneq F_r\in \sd(\Delta)
           \hspace{5pt}|\hspace{5pt}F_r=F,\;r\geq 1\right\}.$$
         By Proposition \ref{Hochster} and Lemma \ref{lemma:no1} it holds that
         \begin{eqnarray*}
            &\empty&\Hilb\;\left ( H^{i} \left ( k[\sd(\Delta)] \right ),t \right )\\	
            &=& \sum_{\sigma\in \sd(\Delta)}\dim_k \widetilde{H}_{i-|\sigma|-1} \left ( 
            \link_{\sd(\Delta)} \left (\sigma \right );k \right ) 
            \left ( \frac{1}{t-1} \right )^{|\sigma|}\\
            &=& \sum_{F \in \Delta} \sum_{r=1}^{|F|} \sum_{\substack{\sigma\in \sd(\Delta)[F] \\
            |\sigma|=r}}\dim_k \widetilde{H}_{i-|\sigma|-1} \left ( \link_{\sd(\Delta)}\left ( 
            \sigma \right );k \right ) \left (\frac{1}{t-1}\right )^{r}\\
            &\empty& + \; \dim_k \widetilde{H}_{i-|\emptyset|-1} \left ( \link_{\sd(\Delta)} 
            \left(\emptyset\right);k\right)\cdot \left (\frac{1}{t-1}\right )^{|\emptyset |}\\
            &=& \sum_{F \in \Delta} \sum_{r=1}^{|F|} \sum_{\substack{\sigma \in \sd(\Delta)[F]\\ 
            |\sigma|=r}} \dim_k \widetilde{H}_{i-|F|-1}\left (\link_{\Delta} F;k \right ) \cdot 
            \left (\frac{1}{t-1} \right )^{r}\\
            &\empty& + \; \dim_k \widetilde{H}_{i-1} \left ( \sd(\Delta);k \right ) \\
            &=& \sum_{F \in \Delta} \sum_{r=1}^{|F|} \Big|\left\{\sigma \in\sd(\Delta)[F]
            \hspace{5pt}|\hspace{5pt}|\sigma|=r\right\}\Big| \cdot \frac{\dim_k \widetilde{H}_{i-|F|-1} 
            \left (\link_{\Delta} F;k\right )}{(t-1)^r}\\
            &\empty& + \; \dim_k \widetilde{H}_{i-1}\left ( \Delta ; k \right )\\
            &=& \sum_{F \in \Delta} \sum_{r=1}^{|F|} f_r^{\sd(\partial F)}
            \cdot \frac{\dim_k 
            \widetilde{H}_{i-|F|-1} \left ( \link_{\Delta} F;k \right )}{(t-1)^r}\\
            &\empty& + \; \dim_k \widetilde{H}_{i-1} \left ( \Delta ; k \right )\\
            &=& \sum_{F \in \Delta} \sum_{r=1}^{|F|} f_{r-2}^{\sd(\partial F)} \cdot 
            \dim_k \widetilde{H}_{i-|F|-1} \left ( \link_{\Delta} F; k\right )\cdot \left 
            (\frac{1}{t-1} \right )^{r}\\
            &\empty & + \; \dim_k \widetilde{H}_{i-1}\left (\Delta ; k \right ).
         \end{eqnarray*}
         Since $f_{r-2}^{\sd(\partial F)}= r!\cdot S\left ( |F|, r \right )$, where 
         $S\left(m,r\right)$ denotes the Stirling number of the second kind (see \cite{S}), 
         it follows that
         \begin{eqnarray*}
            &\empty&\Hilb\;\left ( H^{i} \left ( k[\sd(\Delta)]\right),t \right )\\
            &=&\sum_{F \in \Delta} \sum_{r=1}^{|F|} r! \cdot S\left (|F|,r \right ) 
            \cdot \dim_k \widetilde{H}_{i-|F|-1}\left (\link_{\Delta} F; k \right ) \cdot 
            \frac{1}{(t-1)^r}\\
            &\empty& + \; \dim_k \widetilde{H}_{i-1}\left ( \Delta ; k \right )\\
            &=& \sum_{F \in \Delta} \left ( \sum_{r=1}^{|F|} r! \cdot S\left ( |F|,r \right )
            \cdot \frac{1}{(t-1)^r}\right ) \cdot \dim_k \widetilde{H}_{i-|F|-1}
            \left (\link_{\Delta} F;k \right )\\
            &\empty & + \; \dim_k \widetilde{H}_{i-1} \left ( \Delta ;k \right )\\
%            &=& \sum_{m=1}^{\dim \Delta +1} \sum_{\substack{F \in \Delta \\|F|=m}} 
%            \left (\sum_{r=1}^{m} r! S(m,r)\frac{1}{(t-1)^r}\right ) 
%            \dim_k \widetilde{H}_{i-m-1}\left ( \link_{\Delta}F;k \right )\\
%            &\empty& + \; \dim_k \widetilde{H}_{i-1}\left ( \Delta;k \right )\\
%            &=& \sum_{m=1}^{\dim \Delta +1} \sum_{\substack{F \in \Delta \\ |F|=m}}
%            \left (\sum_{r=1}^{m} r!S\left ( m,r \right )(t-1)^{m-r} \right )
%            \frac{\dim_k \widetilde{H}_{i-m-1} \left (\link_{\Delta} F;k \right )}{(t-1)^m}\\
%            &\empty& + \; \dim_k \widetilde{H}_{i-1} \left ( \Delta ; k \right )\\
            &=&\sum_{m=1}^{\dim \Delta +1} \sum_{\substack{F \in \Delta \\ |F|=m}} 
            \left (\frac{\displaystyle{\sum_{r=1}^{m}}\left (r! S \left (m,r \right ) 
            \sum_{k=0}^{m-r}\binom{m-r}{k}t^{k}(-1)^{m-r-k}\right )}{(t-1)^m}\right ) \cdot \\
            &\empty& ~\hskip6cm \cdot \; \dim_k \widetilde{H}_{i-m-1}\left (\link_{\Delta} F;k\right ) \\
            &\empty& + \; \dim_k \widetilde{H}_{i-1}\left (\Delta; k\right )\\
            &=& \sum_{m=1}^{\dim \Delta +1} \sum_{\substack{F \in \Delta \\ |F|=m}}
            \left (\frac{\displaystyle{\sum_{k=0}^{m-1}}\left (\sum_{r=1}^{m-k}r!S 
            \left ( m,r \right )\binom{m-r}{k} (-1)^{m-r-k} \right ) t^k}{(t-1)^m}\right ) \cdot \\
            &\empty& ~\hskip6cm \cdot \; \dim_k \widetilde{H}_{i-m-1}\left ( \link_{\Delta} F; k\right ) \\
            &\empty& + \; \dim_k \widetilde{H}_{i-1} \left ( \Delta; k \right )\\
            &\begin{array}{c} (*) \\ = \end{array} & \sum_{m=1}^{\dim \Delta +1} 
            \sum_{\substack{F \in \Delta \\ |F|=m}}\left (\frac{\sum_{k=0}^{m-1}
            \left|\left\{\sigma \in S_m \hspace{5pt}|\hspace{5pt} \des(\sigma)=k \right\}\right|
            \cdot t^k}{(t-1)^m}\right ) \cdot \\ 
            &\empty& ~\hskip6cm \cdot \; \dim_k \widetilde{H}_{i-m-1}\left (\link_{\Delta} F;k \right )\\
            &\empty& + \; \dim_k \widetilde{H}_{i-1} \left (\Delta; k \right).
         \end{eqnarray*}
         All manipulations are straight forward, except for (*) which uses a well known
         formula for the Eulerian numbers (see \cite[Corollary 1.18]{Bona}).
      \end{proof}

   \subsection{Depth} \label{Depth}

      \begin{corollary}
         Let $\Delta$ be a simplicial complex. Then
         $$\depth(k[\Delta])= \depth(k[\sd(\Delta)]).$$
      \end{corollary}

      \begin{proof}\label{(1)}
         By a theorem of Grothendieck (see \cite[Theorem 6.2.7]{BS}) the depth of $k[\Delta]$ is 
         given by
         $$\depth(k[\Delta])=\min \left\{i\hspace{5pt}|\hspace{5pt}H^i(k[\Delta])\neq 0\right\}.$$
         By Proposition \ref{Hochster} for the depth of $k[\Delta]$ we get\\
         \begin{eqnarray*}
            \depth(k[\Delta])&=&\min \left\{i\hspace{5pt}|\hspace{5pt}H^i(k[\Delta])\neq 0\right\}\\
            &=&\min \left\{i\hspace{5pt}|\hspace{5pt}\Hilb\;(H^i(k[\Delta]),t)\neq 0\right\}\\
            &=&\min \left\{i\hspace{5pt}|\hspace{5pt}\exists F \in \Delta : \dim_k \widetilde{H}_{i-
            \left|F\right|-1}\left(\link_{\Delta}F;k\right)\neq 0 \right\}.\\
         \end{eqnarray*}
         Analogously,
         $\depth(k[sd(\Delta)])=$ $$\min \left\{i\hspace{5pt}|\hspace{5pt}\exists 
         \sigma\in \sd(\Delta) : \dim_k\widetilde{H}_{i-|\sigma|-1}\left(lk_{\sd(\Delta)}(\sigma);k\right)
         \neq 0\right\}.$$
         By Lemma \ref{lemma:no1} for $\sigma = F_1 \subsetneq \cdots \subsetneq F_r := F$
         we have

         \begin{eqnarray*}
            &\empty&\widetilde{H}_{m-|\sigma|-1}\left(\link_{\sd(\Delta)}\left(\sigma \right);k\right) = 
            \widetilde{H}_{m-|F|-1}\left(\link_{\Delta}F;k\right),
         \end{eqnarray*}
         which implies the assertion.
      \end{proof}

   \subsection{Projective Dimension} \label{Projective-Dimension}

     We denote by $\prdim(k[\Delta])$ the projective dimension of $k[\Delta]$. 
     For a simplicial complex $\Delta$ over ground set $[n]$ and $f_0^\Delta = n$ we
     obtain using The Auslander-Buchsbaum formula (see \cite[Theorem 19.9]{E}):
     \begin{eqnarray*} 
        \prdim(k[\Delta]) & = & f_0^\Delta - \depth(k[\Delta]) \\
            & = & \sum_{i \geq 0} f_i^\Delta - \sum_{i \geq 1} f_i^\Delta - \depth(k[\sd(\Delta)])\\
            & = & f_0^{\sd(\Delta)} - \depth(k[\sd(\Delta)]) - \sum_{i \geq 1} f_i^\Delta \\
            & = & \prdim(k[\sd(\Delta)]) - \sum_{i \geq 1} f_i^\Delta.
     \end{eqnarray*}

   \subsection{Regularity} \label{Regularity}

      \begin{proposition}
         Let $\Delta$ be a simplicial complex on vertex set $[n]$.
         \begin{eqnarray*}
            \reg (k[\Delta]) \leq \reg (k[\sd(\Delta)])=
            \left\{ 
            \begin{array}{lcl}
              \dim \Delta \hspace{3pt} & \text{if } & \widetilde{H}_{\dim \Delta}\left(\Delta;k\right)=0\\
              \dim \Delta +1 \hspace{1pt} & \text{if } & \widetilde{H}_{\dim \Delta}\left(\Delta;k\right) 
              \neq 0\\
            \end{array} 
            \right.
         \end{eqnarray*}
         Moreover, if $\widetilde{H}_{\dim \Delta}\left(\Delta;k\right) \neq 0$ then 
         $$\reg (k[\Delta])= \reg (k[\sd(\Delta)]) = \dim \Delta +1.$$
      \end{proposition}

      \begin{proof}
        We use the following characterization of regularity \cite{BS}

        $$\reg(k[\Delta]) = \max\{i+j~|~H^i(k[\Delta])_j \neq 0\}.$$

        \noindent By Proposition \ref{Hochster} we have that
        \begin{eqnarray*}
           \Hilb \left(H^i\left( k[\Delta] \right),t \right)&=& \sum_{l \in \integers}
           \dim_k H^i \left(k[\Delta] \right)_l \cdot t^l\\
           &=& \sum_{F \in \Delta} \dim_k \widetilde{H}_{i-|F|-1}\left(\link_{\Delta}F;k \right) 
           \cdot \left(\frac{1}{t-1}\right)^{|F|}\\
           &=& \sum_{F \in \Delta} \dim_k \widetilde{H}_{i-|F|-1}\left(\link_{\Delta}F;k\right)\cdot 
           \left( \sum_{n \geq 1} \left( \frac{1}{t} \right)^n \right)^{|F|}\\
           &=& \sum_{F \in \Delta} \dim_k \widetilde{H}_{i-|F|-1} \left(\link_{\Delta}F;k \right) 
           \cdot \sum_{a \in \left( \mathbb{Z}_{-}\setminus \left\{0\right\}\right)^{|F|}} t^{|a|}. \\
        \end{eqnarray*}
        Here $|a|=a_1+ \ldots +a_{|F|}$ and $\integers_{-}=\left\{0,-1,-2,\ldots\right\}$.
        We conclude
        \begin{eqnarray*}
           &\empty& H^i\left(k[\Delta]\right)_l \neq 0 \\
           &\Longleftrightarrow &\dim_k H^i \left(k[\Delta]\right)_l \neq 0\\
           &\Longleftrightarrow& \exists a \in \mathbb{Z}_{-}^{n},\; |a|=l \text{ and } 
           \exists \hspace{5pt} F \in \Delta \text{ such that } \Big|\left\{s\hspace{5pt}|\hspace{5pt} 
           a_s<0\right\}\Big|=|F| \\
           &\empty&\text{ and }\dim_k \widetilde{H}_{i-|F|-1}\left( \link_{\Delta} F;k \right) \neq 0
        \end{eqnarray*}
        As usual, for a $\integers$-graded module $M = \bigoplus_{m \in \integers} M_m$, we write 
        $\End(M)$ for $\sup\left\{m\in \integers\hspace{5pt}|\hspace{5pt}M_m\neq 0\right\}$.
        The above directly yields $\End( H^i\left(k[\Delta]\right))\leq 0$ and choosing 
        $a \in \left\{0,-1\right\}^n$ we see that 
        $\End (H^i\left(k[\Delta]\right))\geq -(\dim \Delta +1)$ if 
        $H^i\left(k[\Delta]\right)\neq 0$.\\
        More precisely, if $H^i\left(k[\Delta]\right)\neq 0$ then $\End( H^i(k[\Delta]))$ is given by
        \begin{eqnarray*}
           \End (H^i\left(k[\Delta]\right)) &=& \sup \left\{m \in \integers
           \hspace{5pt}|\hspace{5pt}H^i\left(k[\Delta]\right)_m \neq 0 \right\}\\
           &=&\sup \left\{0\geq m \geq -(\dim \Delta +1)\hspace{5pt}|\hspace{5pt}H^i\left(k[\Delta]
           \right)_m \neq 0\right\}\\
           &=&\sup \left\{0 \geq m \geq -(\dim \Delta+1) \hspace{5pt}\big|\substack{ 
           \exists F \in \Delta,\;|F|=-m:\\ \widetilde{H}_{i-|F|-1}\left(\link_{\Delta}F;k\right) 
           \neq 0}\right\}\\
           &=& \sup \left\{-|F|\hspace{5pt}|\hspace{5pt}F \in \Delta:\;\widetilde{H}_{i-|F|-1}
           \left(\link_{\Delta}F;k\right)\neq 0\right\}\\
           &=&- \inf \left\{|F|\hspace{5pt}|\hspace{5pt}F\in \Delta: \; 
           \widetilde{H}_{i-|F|-1}\left(\link_{\Delta}F;k\right)\neq 0\right\}.
        \end{eqnarray*}
        We now show $\reg (k[\Delta]) \leq \reg (k[\sd(\Delta)])$.\\
        From the previous consideration and $|F_1\subsetneq \ldots \subsetneq F_t|\leq|F_t|$ 
        it follows that 
        \begin{eqnarray*}
           \End (H^i\left(k[\sd(\Delta)]\right))&=& 
           - \inf \left\{~t~\big|\substack{\exists F_1\subsetneq \ldots \subsetneq F_t \in 
           \sd(\Delta) :\\ \widetilde{H}_{i-t-1}\left(\link_{\sd(\Delta)}\left(F_1\subsetneq \ldots 
           \subsetneq F_t\right);k\right)\neq 0}\right\}\\
           &=&-\inf \left\{~t~\big|\substack{\exists F_1\subsetneq \ldots \subsetneq 
           F_t\in \sd(\Delta):\\ \widetilde{H}_{i-|F_t|-1}\left(\link_{\Delta}F_t;k\right)\neq 0}\right\}\\
           &\geq&-\inf\left\{|F_t|\big|\substack{\exists F_1\subsetneq \ldots \subsetneq F_t \in 
           \sd(\Delta):\\ \widetilde{H}_{i-|F_t|-1}\left(\link_{\Delta}F_t;k\right)\neq 0}\right\}\\
           &=& - \inf \left\{~|F_t|~\big|F_t \in \Delta:\widetilde{H}_{i-|F_t|-1}
           (\link_{\Delta}F_t;k)\neq 0\right\}\\
           &=& \End (H^i\left(k[\Delta]\right)).
        \end{eqnarray*}
        Thus $ \End (H^i\left(k[\Delta]\right))+i\leq \End (H^i\left(k[\sd(\Delta)]\right))+i$ 
        for all $i$. It follows that $\reg (k[\Delta])\leq \reg (k[\sd(\Delta)])$, as desired.\\
        We claim $\reg (k[\sd(\Delta)])= \left\{ \begin{array}{lcl}
        \dim \Delta & \text{if } & \widetilde{H}_{\dim \Delta}\left(\Delta;k\right)=0\\
        \dim \Delta+1 & \text{if } & \widetilde{H}_{\dim\Delta}\left(\Delta;k\right)\neq 0\\
        \end{array} \right.$\\
        {\sf Case 1:} $\widetilde{H}_{\dim \Delta}\left(\Delta;k\right)=0.$\\
            By Grothendieck's Non-Vanishing Theorem (see \cite[6.1.4]{BS}) we have that 
            $ H^{\dim k[\Delta]}\left(k[\Delta]\right)=H^{\dim \Delta +1}\left(k[\Delta]\right)\neq 0$.
            Along with the above consideration we conclude that there exists a face $F \in \Delta$ 
            such that $\widetilde{H}_{\dim \Delta +1-|F|-1}\left(\link_{\Delta}F;k\right)=
            \widetilde{H}_{\dim \Delta-|F|}\left(\link_{\Delta}F;k\right)\neq 0$. Since 
            $\widetilde{H}_{\dim \Delta}\left(\Delta;k\right)=0$ this face cannot be the empty face.\\
            Therefore, we can consider $F$ as a one-element flag in $\sd(\Delta)$. 
            From Lemma \ref{lemma:no1} we deduce 
            $$\widetilde{H}_{(\dim \Delta+1)-\underbrace{|F|}_{= 1 \text{~in~}\sd(\Delta)}-1}
               \left(\link_{\sd(\Delta)}F;k\right)
               =\widetilde{H}_{\dim \Delta-\underbrace{|F|}_{\text{in }\Delta}}
               \left(\link_{\Delta}F;k\right)\neq 0. $$
            Thus 
            \begin{eqnarray}
               &\empty& \End (H^{\dim \Delta +1}\left(k[\sd(\Delta)]\right))\nonumber\\
               &=& - \inf \left\{|F_1 \subsetneq \ldots \subsetneq F_t|\big|\substack{F_1 
               \subsetneq \ldots \subsetneq F_t \in \sd(\Delta):\\  \widetilde{H}_{\dim \Delta -
               |F_1 \subsetneq \ldots \subsetneq F_t|}\left(\link_{\sd(\Delta)}
               \left(F_1\subsetneq \ldots \subsetneq F_t\right);k\right)\neq 0}\right\}\nonumber\\
               &\geq& -\underbrace{|F|}_{= 1 \text{~in~}\sd(\Delta)}=-1. \label{eqn:no1}
            \end{eqnarray}
            This implies 
            \begin{eqnarray*}
               \reg (k[\sd(\Delta)])&\geq& \End (H^{\dim \Delta +1}\left(k[\sd(\Delta)]\right))+
               \dim \Delta +1\\
               &\geq& -1+\dim \Delta +1 = \dim \Delta .\\
            \end{eqnarray*}
            We also know that $\End (H^i\left(k[\sd(\Delta)]\right))\leq 0$. Hence
            $$\End (H^i\left(k[\sd(\Delta)]\right))+i\leq i < \dim \Delta +1,\; 0\leq i \leq 
            \dim \Delta.$$
            From Inequality (\ref{eqn:no1}) we deduce 
            $0\geq \End (H^{\dim \Delta +1}\left(k[\sd(\Delta)]\right))\geq -1$.
            Since
            $$\widetilde{H}_{\dim \Delta-|\emptyset|}\left(\link_{\sd(\Delta)}(\emptyset);k\right)=
               \widetilde{H}_{\dim \Delta}\left(\sd(\Delta);k\right)=
               \widetilde{H}_{\dim \Delta}\left( \Delta;k\right)=0$$
            we conclude that $\End (H^{\dim \Delta +1}\left(k[\sd(\Delta)]\right))\neq 0$.
            Therefore, $$\End (H^{\dim \Delta +1}\left(k[\sd(\Delta)]\right))=-1.$$
            Thus $\End (H^{\dim \Delta +1}\left(k[\sd(\Delta)]\right))+\dim \Delta +1=\dim \Delta$.
            This finally proves the claim for this case.\\
        {\sf Case 2:} 
            $\widetilde{H}_{\dim \Delta}\left(\Delta;k\right)\neq 0$.\\
            From the fact that $\End(H^i(k[\Delta]) \leq 0$ it follows that
            $\End (H^i\left(k[\Delta]\right))+i\leq \dim \Delta$, $0\leq i \leq \dim \Delta$.\\
            By $\widetilde{H}_{\dim \Delta}\left(\Delta;k\right)\neq 0$ we get
            \begin{eqnarray*}
               0 & \geq & \End (H^{\dim \Delta +1}\left(k[\Delta]\right)) \\
               \empty & = & -\inf \left\{|F|~\big|~ F \in \Delta:\; 
               \widetilde{H}_{\dim \Delta +1-|F|}\left(\link_{\Delta}\left(F\right);k\right)\neq 0
               \right\}
               \geq |\emptyset|=0.
           \end{eqnarray*}
            Therefore, 
            \begin{eqnarray*}
               \reg (k[\Delta])&=&\sup \left\{\End (H^i(k[\Delta]))+i\hspace{5pt}|
               \hspace{5pt}0 \leq i \leq \dim \Delta +1\right\}\\
               &=& \sup \Big(\left\{\underbrace{\End (H^i\left(k[\Delta]\right))+i}_{\leq \dim \Delta}                \hspace{5pt} | \hspace{5pt} 0\leq i \leq \dim \Delta \right\}\\
               &\empty&\cup \left\{\underbrace{\End (H^{\dim \Delta +1}\left(k[\Delta]\right))+ 
               \dim \Delta +1}_{= \dim \Delta +1}\right\}\Big)\\
               &=& \dim \Delta +1.
            \end{eqnarray*}
            By $\dim \Delta = \dim \sd(\Delta)$ and $\widetilde{H}_{\dim \Delta}(\Delta) = 
            \widetilde{H}_{\dim \Delta}(\sd(\Delta))$ the above also applies to $\sd(\Delta)$ and 
            $\reg (k[\sd(\Delta)])= \dim \Delta +1$ follows.
      \end{proof}

     \subsection{Height and Multiplicity} \label{Height-Multiplicity}

        \begin{proposition} \label{height, multiplicity}
           Let $\Delta$ be a simplicial complex with $f$-vector 
           $\fvec^\Delta=\left(f_0^{\Delta},\ldots,f_{\dim \Delta}^{\Delta}\right)$.
           Then
           \begin{eqnarray*}
              \height (I_{\sd(\Delta)})&=&\sum_{l=0}^{\dim \Delta}(f_l^{\Delta}-1)\text{   and}\\
              e(k[\sd(\Delta)])&=&(\dim \Delta +1)!\cdot f_{\dim \Delta}^{\Delta}.
           \end{eqnarray*}
         \end{proposition}

         \begin{proof}
            From a result by Herzog \& Srinivasan (\cite[2. Antichains]{SH}) one deduces that 
            $\height (I_{\Delta})=f_0^{\Delta}-(\dim \Delta +1)$ and 
            $e(k[\Delta])=f_{\dim \Delta}^{\Delta}$. Thus
            \begin{eqnarray*}
               \height (I_{\sd(\Delta)})&=&f_0^{\sd(\Delta)}-(\dim \sd(\Delta) +1)\\
               &=&\sum_{l=0}^{\dim \Delta}f_l^{\Delta}-\dim \Delta-1\\
               &=& \sum_{l=0}^{\dim \Delta}(f_l^{\Delta}-1).
            \end{eqnarray*}
            A simple counting argument shows that 
            \begin{eqnarray*}
               e(k[\sd(\Delta)])&=&f_{\dim \Delta}^{\sd(\Delta)}\\
               &=&(\dim \Delta +1)!\cdot f_{\dim \Delta}^{\Delta}.
            \end{eqnarray*}

         \end{proof}
  \subsection{Cohen-Macaulay-ness} \label{Cohen-Macaulay-ness}
     A simplicial complex is Cohen-Macaulay over a field $k$ if $k[\Delta]$ is a
     Cohen-Macaulay ring (see \cite[Chapter 5]{BH} for background on Cohen-Macaulay
     simplicial complexes).
     It is a well known fact from geometric combinatorics that Cohen-Macaulay-ness over a
     field $k$ of a simplicial complex depends on its topological realization only
     (see \cite{B}). Since
     $\Delta$ and $\sd(\Delta)$ have homeomorphic geometric realizations it follows that
     $\Delta$ is Cohen-Macaulay over $k$ if and only if $\sd(\Delta)$ is. 
 
  \subsection{Koszulness} \label{Koszulness}
     The minimal nonfaces of $\sd(\Delta)$ are of cardinality two -- the pairs of
     faces of $\Delta$ that are incomparable. Therefore, the Stanley-Reisner ideal
     of $\sd(\Delta)$ is generated by (squarefree) monomials of degree $2$ and hence by
     a result of Fr\"oberg (see \cite{F} for the result and background on Koszul algebras)
     it is Koszul.

  \subsection{Golod-ness} \label{Golod-ness}
     We have already seen in Section \ref{Koszulness} that $I_\Delta$ is generated by
     squarefree monomials of degree two. By a result of Berglund \& J\"ollenbeck 
     \cite[Theorem 7.4]{BJ} we know
     that in this situation $k[\sd(\Delta)]$ is Golod if and only if the $1$-skeleton of
     $\sd(\Delta)$ is a chordal graph; i.e. any cycle of length $\geq 4$ has chord.
     We refer the reader to \cite{GL} for background on Golod-ness.
     Now assume the $1$-skeleton of $\Delta$ has a chordless cycle of length $\ell \geq 3$ 
     -- here we regard triangles as chordless cycles.
     Then after barycentric subdivision this chordless cycle turns into a chordless cycle
     of length $2\ell \geq 6$. Hence, $\sd(\Delta)$ cannot be Golod. So assume the 
     $1$-skeleton of $\Delta$ has no chordless cycle of length $\geq 3$. Then $\dim \Delta \leq 1$
     and $\Delta$ is a graph. Having no chordless cycle of length $\geq 3$ then implies that $\Delta$ has no
     cycle and hence is a forest. Now the barycentric subdivision of a forest is a forest and hence has
     no cycles which implies that $\sd(\Delta)$ is Golod. Thus:

     \begin{proposition}
       Let $\Delta$ be a simplicial complex. Then 
       $\sd(\Delta)$ is Golod if and only if $\Delta$ is a forest.
    \end{proposition}

\section{Auxiliary Lemmas and Inequalities} \label{AuxiliaryLemmas}

  The basic result which allows us to verify the Multiplicity Conjecture for 
  barycentric subdivisions is the following classical theorem by Hochster which
  expresses the Betti numbers of $k[\Delta]$ in terms of homology groups of
  restrictions of $\Delta$. For a simplicial complex on ground set $\Omega$ the
  restriction $\Delta_W$ of $\Delta$ to a subset $W \subseteq \Omega$ is 
  the simplicial complex $\Delta_{W}:=\left\{ F \in\Delta\hspace{5pt}|
  \hspace{5pt} F \subseteq W \right\}$.

  \begin{proposition}[\cite{H}]\label{Hochster_2}
    Let $\Delta$ be a simplicial complex on vertex set $[n]$ and
    let $\beta_{ij}$ be the bigraded Betti-numbers of $k[\Delta]$.   
%    Let $\Hilb\;(T_i,\lambda)=\sum_{a\in \integers^{n}}\beta_{ia}\cdot \lambda^a$ the fine 
%    Hilbert series of $T_i=Tor_i^T\left(k,k[\Delta]\right)$ where $T=k[x_1,\ldots,x_n]$. 
%    Then 
%    \begin{eqnarray*}
%      \Hilb\;(T_i,\lambda)&=&\sum_{\alpha \in \integers^{n}}\dim_k Tor_i^T
%      \left(k[\Delta];k\right)_{\alpha}\cdot \lambda^{\alpha}\\
%      &=&\sum_{\alpha \in \integers^{n}}\beta_{ia}\cdot \lambda^a\\
%      &=&\sum_{W\subseteq [n]}\dim_k \widetilde{H}_{|W|-i-1}\left(\Delta_W;k\right)\cdot 
%      \prod_{i\in W}\lambda_{i}.
%    \end{eqnarray*}
%    In particular, for $i,j\in \mathbb{N}$ 
    Then for $i,j\in \mathbb{N}$ 
    $$\beta_{ij}=\sum_{\substack{W\subseteq [n]\\|W|=j}}
    \dim_k\widetilde{H}_{|W|-i-1}\left(\Delta_W;k\right).$$
  \end{proposition}

  The following corollary is a direct consequence of Proposition \ref{Hochster_2}.

  \begin{corollary}\label{Betti_0}
     Let $\Delta$ be a simplicial complex on vertex set $[n]$ and $\beta_{ij}$ the 
     bigraded Betti numbers of $k[\Delta]$. Then it holds that
     $$\beta_{ij}\neq 0 \hspace{5pt} \Longleftrightarrow \hspace{5pt}
     \exists W\subseteq [n],\;|W|=j\text{ such that }\widetilde{H}_{j-i-1}\left(\Delta_W;k\right)\neq 0.$$
  \end{corollary}

  \begin{lemma} \label{lemma: Betti_1}
     Let $\Delta$ be a simplicial complex and $\beta_{ij}$ the bigraded Betti numbers of 
     $k[\sd(\Delta)]$.
     Then
     \begin{enumerate}[(i)]
        \item For $1<m< \dim \Delta$ and $2^{m+1}-2-m\leq i<2^{m+2}-2-(m+1)$ we have 
              $\beta_{i,i+m}\neq 0$. 
        \item If $1<\dim \Delta$ then $\prdim(k[\sd(\Delta)]) \geq 4$. Equivalently,  
              for $1\leq i\leq 4$ there exist $k_i\geq 1$ such that 
              $\beta_{i,i+k_i}\neq 0$.
     \end{enumerate}
  \end{lemma}

  \begin{proof}
%     First, we prove (ii). 
%     Since $I_{\sd(\Delta)}$ is generated by monomials of degree two we have 
%     $\beta_{1,2}\neq 0$. If $\sum_{j\geq i}\beta_{i,j}\neq 0$ for $i\geq2$ there exist 
%     $k_i\geq 1$ with $\beta_{i,i+k_i}\neq 0$. \\
     Clearly, a simplicial complex $\Delta$ contains at least one $(\dim \Delta)$-simplex. Thus 
     $\sum_{l=0}^{\dim \Delta}f_l^{\Delta}\geq 2^{\dim \Delta+1}-1$. 
     From Section \ref{Projective-Dimension} and by $\dim k[\Delta] = \dim (\Delta)+1$ we deduce
     \begin{eqnarray*}
        \prdim(k[\sd(\Delta)])&=& \prdim(k[\Delta]) + \sum_{i \geq 1} f_i^\Delta \\
        & = & f_0^{\Delta}-\depth (k[\Delta]) + \sum_{i \geq 1} f_i^\Delta \\
        &\geq& \sum_{l=0}^{\dim \Delta} f_l^{\Delta}-\dim \Delta-1 \\
     \end{eqnarray*}
     In particular, if $\dim(\Delta) > 1$ then $\prdim(k[\sd(\Delta)]) \geq 4$, which proves (ii). 

     The same inequality shows that if $\dim \Delta \geq m+1$ then
     $\prdim(k[\sd(\Delta)])$ $\geq 2^{m+2}-2-(m+1)$.
     Thus for $1<m<\dim\Delta$ and for $2^{m+1}-2-m\leq i<2^{m+2}-2-(m+1)$ there 
     exist $j\geq 1$ such that $\beta_{i,i+j}\neq 0$.\\

     Let $F\in \Delta$ with $\dim F=m$. The boundary 
     $\partial F$ of $F$ is homeomorphic to an $(m-1)$-sphere. 
     Hence $\sd(\partial F)\cong \sd(\Delta)_{\mathring{\partial F}}\cong \mathbb{S}^{m-1}$. 
     Thus
     $$k=\widetilde{H}_{m-1}\left(\sd(\Delta)_{\mathring{\partial F}};k\right)=
     \widetilde{H}_{|\mathring{\partial F}|-(|\mathring{\partial F}|-m)-1}\left(\sd(\Delta)_{\mathring{\partial F}};k\right).$$
     Using Corollary \ref{Betti_0} we conclude
     $$\beta_{|\mathring{\partial F}|,|\mathring{\partial F}|-m}\stackrel{|\mathring{\partial F}|=2^{m+1}-2}{=
     }\beta_{2^{m+1}-2-m,2^{m+1}-2}\neq 0.$$
     It now remains to show that $\beta_{i,i+m}\neq 0$ for $1<m<\dim \Delta$ and $2^{m+1}-2-m<i<2^{m+2}-2-(m+1)$. 
     Since $\dim \Delta \geq m+1$ there exists an $(m+1)$-dimensional face $G\in \Delta$. 
     Choose $v \in G$ and set $F = G \setminus \{v\}$. Then
     $\sd(\Delta)_{\mathring{\partial F}}\cong \mathbb{S}^{m-1}$ and
     $\sd(\Delta)_{\mathring{\partial F}\cup H}\simeq \mathbb{S}^{m-1}$ for all 
     $H\subseteq\displaystyle{\mathring{\Delta}_G}\setminus \left\{G,F,v\right\}$ and $H\cap \partial F=\emptyset$. 
     The last assertion holds 
     because when restricting $\sd(\Delta)$ to $\mathring{\partial F}\cup H$ we are only 
     adding simplices to $\sd(\Delta)_{\mathring{\partial F}}$ which can be contracted to 
     $\sd(\Delta)_{\mathring{\partial F}}$. Moreover, since we do not add $v$, $F$ and $G$ the complex 
     $\sd(\Delta)_{\mathring{\partial F}\cup H}$ still contains the cycle induced by $\partial F$. We conclude
     \begin{eqnarray*}
        \widetilde{H}_{m-1}\left(\sd(\Delta)_{\mathring{\partial F}\cup H};k\right)&=&\widetilde{H}_{|H|+
        |\mathring{\partial F}|-(|H|+|\mathring{\partial F
        }|-m)-1}\left(\sd(\Delta)_{\mathring{\partial F}\cup H};k\right)\\
        &=&\widetilde{H}_{|H|+|\mathring{\partial F}|-(|H|+|\mathring{\partial 
        F}|-m)-1}\left(\mathbb{S}^{m-1};k\right)\\
        &=&\widetilde{H}_{m-1}\left(\mathbb{S}^{m-1};k\right)\neq 0.
     \end{eqnarray*}
     By Corollary \ref{Betti_0} we obtain $\beta_{|H|+|\mathring{\partial F}|-m,|H|+|\mathring{\partial F}|}\neq 0$.\\
     From $$2^{m+1}-2\stackrel{H=\emptyset}{\leq}|H|+|\mathring{\partial F}|
     \stackrel{H\cup \mathring{\partial F}=\mathring{\Delta}_G\setminus \left\{G,F,v\right\}}{\leq}2^{m+2}-4$$ 
     it follows that $\beta_{j-m,j}\neq 0$ for $2^{m+1}-2\leq j \leq 2^{m+2}-4$.\\
     Thus $\beta_{i,i+m}\neq 0$ for $2^{m+1}-2-m\leq i \leq 2^{m+2}-4-m<2^{m+2}-2-(m+1)$, as desired.
  \end{proof}

  \begin{lemma} \label{lemma: Betti_2}
     Let $\Delta$ be a simplicial complex such that $\widetilde{H}_{\dim \Delta}\left(\Delta;k\right)=0$ and 
     let $\beta_{ij}$ be the bigraded Betti numbers of $k[\sd(\Delta)]$.
     Then 
     $$\beta_{i,i+\dim \Delta}\neq 0\hspace{5pt} \mbox{~for~} 
     2^{\dim \Delta+1}-2-\dim \Delta \leq i \leq \sum_{j=0}^{\dim \Delta}(f_j^{\Delta}-1).$$
  \end{lemma}

  \begin{proof}\label{Betti}
     Let $F \in \Delta$ with $\dim F =\dim \Delta$. We show that 
     $\widetilde{H}_{\dim \Delta -1}\left(\sd(\Delta)_{\mathring{\Delta}\setminus
     \left\{F\right\}};k\right)\neq 0$.
     By an elementary homotopy $\sd(\Delta)_{\mathring{\Delta}\setminus \left\{F\right\}} \simeq 
     \sd(\Delta)\setminus \left\{F\right\}$. Now 
     consider the long exact sequence in homology of the pair 
     $\left(\sd(\Delta),\sd(\Delta)\setminus \left\{F\right\}\right)$.
     $$\rightarrow \underbrace{\widetilde{H}_{\dim \Delta}\left(\sd(\Delta);k\right)}_{=0 \text { by assumption}}\stackrel{q_{\dim \Delta}}{\rightarrow}
     \widetilde{H}_{\dim \Delta}\left(\sd(\Delta),\sd(\Delta)\setminus \left\{F\right\};k\right)\stackrel{\partial}{\rightarrow} $$
     $$\widetilde{H}_{\dim \Delta-1}\left(\sd(\Delta) \setminus \left\{F\right\};k\right)\stackrel{i}{\rightarrow}
     \widetilde{H}_{\dim \Delta -1}\left(\sd(\Delta);k\right)\stackrel{q_{\dim \Delta -1}}{\rightarrow}$$
     $$\widetilde{H}_{\dim \Delta -1}\left(\sd(\Delta),\sd(\Delta)\setminus \left\{F\right\};k\right)\rightarrow \ldots$$
     Since $\left(\sd(\Delta),\sd(\Delta)\setminus F\right)$ is a good pair we have
     $$\widetilde{H}_{\dim \Delta}\left(\sd(\Delta),\sd(\Delta)\setminus \left\{F\right\};k\right) \cong 
     \widetilde{H}_{\dim \Delta}\left(\sd(\Delta)/(\sd(\Delta)\setminus\left\{F\right\});k\right)
     $$
     $$\cong\widetilde{H}_{\dim \Delta}\left(\mathbb{S}^{\dim \Delta};k\right) = k.$$
     The same argument shows $\widetilde{H}_{\dim \Delta -1}\left(\sd(\Delta),\sd(\Delta)\setminus \left\{F\right\};k\right)=0$.
     Along with the above sequence being exact this implies 
     $$\IM\;i=\KER\;q_{\dim \Delta -1}=\widetilde{H}_{\dim \Delta -1}\left(\sd(\Delta);k\right).$$
     Thus
     \begin{eqnarray*}
        \widetilde{H}_{\dim \Delta -1}\left(\sd(\Delta);k\right)&\cong&\widetilde{H}_{\dim \Delta -1}
        \left(\sd(\Delta)\setminus \ \left\{F\right\};k\right)/\KER\; i\\
        &=&\widetilde{H}_{\dim \Delta -1}\left(\sd(\Delta)\setminus \left\{F\right\};k\right)/\IM  \;\partial.
     \end{eqnarray*}
     Since the above sequence is exact it holds that $0=\IM \; q_{\dim \Delta}=
    \KER \; \partial$. This yields 
     $$\IM \;\partial\cong \widetilde{H}_{\dim \Delta}\left(\sd(\Delta),\sd(\Delta)\setminus \left\{F\right\};k\right)\cong k.$$  
     Thus 
     \begin{eqnarray*}
        \widetilde{H}_{\dim \Delta -1}\left(\sd(\Delta)\setminus \left\{F\right\};k\right)&\cong&\widetilde{H}_{\dim \Delta -1}
        \left(\sd(\Delta);k\right)\oplus \IM\;\partial\\
        &\cong&\widetilde{H}_{\dim \Delta -1}\left(\sd(\Delta);k\right)\oplus k\neq 0.
     \end{eqnarray*}
     It follows that $\beta_{\sum_{j=0}^{\dim \Delta}(f_j^{\Delta}-1),\sum_{j=0}^{\dim \Delta}f_j^{\Delta}-1}=
     \beta_{|\mathring{\Delta}\setminus \left\{F\right\}|-\dim \Delta,|\mathring{\Delta}\setminus \left\{F\right\}|}\neq 0$.
     Our next aim is to prove that $\widetilde{H}_{\dim \Delta -1}\left(\sd(\Delta)_{\mathring{\Delta}\setminus (A\cup\left\{F\right\})};k\right)
     \neq 0$ for  $A\subseteq\Delta\setminus \Delta_F$. First using induction on the cardinality of $A$ we show, that 
     $\widetilde{H}_{\dim \Delta}\left(\sd(\Delta)_{\mathring{\Delta}\setminus A};k\right)=0$.\\
     For $|A|=0$ this is our assumption on the homology of $\Delta$.
     Assume $|A| \geq 1$.  Let $\widetilde{H}_{\dim \Delta}\left(\sd(\Delta)_{\mathring{\Delta}\setminus A};k\right)=0$ for all $A\subseteq \Delta\setminus \Delta_F$ with $|A|=m$ and
     let $B\subseteq \Delta\setminus \Delta_F$ with $|B|=m+1$. Consider $A:=B\setminus \left\{v\right\}$ for some $v\in B$. 
     By the induction hypothesis we have $\widetilde{H}_{\dim \Delta}\left(\sd(\Delta)_{\mathring{\Delta}\setminus A};k\right)=0$.
     Consider the exact homology sequence of the $\left(\sd(\Delta)_{\mathring{\Delta}\setminus A},\sd(\Delta)_{\mathring{\Delta}\setminus (A\cup \left\{F\right\})}\right)$: 
     \begin{eqnarray*}
        \cdots\widetilde{H}_{\dim \Delta +1}\left(\sd(\Delta)_{\mathring{\Delta}\setminus A},\sd(\Delta)_{\mathring{\Delta}
        \setminus (A\cup \left\{v\right\})};k\right)\stackrel{\partial}{\rightarrow}\widetilde{H}_{\dim \Delta}
        \left(\sd(\Delta)_{\mathring{\Delta}\setminus (A\cup \left\{v\right\})};k\right)\\
        \stackrel{i}{\rightarrow} \underbrace{\widetilde{H}_{\dim \Delta}
        \left(\sd(\Delta)_{\mathring{\Delta}\setminus A};k\right)}_{=0 \text{ by induction hypothesis}}\rightarrow
        \widetilde{H}_{\dim \Delta}\left(\sd(\Delta)_{\mathring{\Delta}\setminus A},\sd(\Delta)_{\mathring{\Delta}
        \setminus (A\cup\left\{v\right\})};k \right)\cdots 
     \end{eqnarray*}
     Since $\sd(\Delta)_{\mathring{\Delta}\setminus A}$ and $\sd(\Delta)_{\mathring{\Delta}\setminus (A\cup \left\{v\right\})}$ 
     are $(\dim \Delta)$-dimensional CW-com\-ple\-xes it follows that $\sd(\Delta)_{\mathring{\Delta}\setminus A}/\sd(\Delta)_{\mathring{\Delta}
     \setminus (A\cup \left\{v\right\})}$ is a $(\dim \Delta)$-dimensional CW-complex. In particular, the complex 
     has no cells in dimension $\dim \Delta +1$. Thus 
     \begin{eqnarray*}
        &\empty& \widetilde{H}_{\dim \Delta +1}\left(\sd(\Delta)_{\mathring{\Delta}\setminus A},
        \sd(\Delta)_{\mathring{\Delta}\setminus (A\cup \left\{v\right\})};k\right) \\ 
      &=&  \widetilde{H}_{\dim \Delta +1}\left(\sd(\Delta)_{\mathring{\Delta}\setminus A}/\sd(\Delta)_{\mathring{\Delta}
        \setminus (A\cup \left\{v\right\})};k\right)=0.
    \end{eqnarray*}
     Exactness of the preceding sequence implies the desired fact 
$$        \widetilde{H}_{\dim \Delta}\left(\sd(\Delta)_{\mathring{\Delta}\setminus (A\cup \left\{v\right\})};k\right)
        =\widetilde{H}_{\dim \Delta}\left(\sd(\Delta)_{\mathring{\Delta}\setminus B};k\right)=0. $$
     Consider the long exact sequence of the pair $\left(\sd(\Delta)_{\mathring{\Delta}\setminus A},
     \sd(\Delta)_{\mathring{\Delta}\setminus (A\cup \left\{F\right\})}\right)$ for an arbitrary $A\subseteq\Delta\setminus \Delta_F$.
     \begin{eqnarray*}
        &\ldots& \widetilde{H}_{\dim \Delta}\left(\sd(\Delta)_{\mathring{\Delta}\setminus A};k\right)\rightarrow
        \widetilde{H}_{\dim\Delta}\left(\sd(\Delta)_{\mathring{\Delta}\setminus A},\sd(\Delta)_{\mathring{\Delta}
        \setminus (A\cup\left\{F\right\})};k\right)\\
        &\rightarrow&\widetilde{H}_{\dim \Delta -1}\left(\sd(\Delta)_{\mathring{\Delta}\setminus (A\cup\left\{F\right\})};k\right)\rightarrow
        \widetilde{H}_{\dim \Delta -1}\left(\sd(\Delta)_{\mathring{\Delta}\setminus A};k\right)\\
        &\rightarrow&\widetilde{H}_{\dim \Delta-1}\left(\sd(\Delta)_{\mathring{\Delta}\setminus A},
        \sd(\Delta)_{\mathring{\Delta}\setminus (A\cup \left\{F\right\})};k\right)\rightarrow \ldots
     \end{eqnarray*}
     Since $\left(\sd(\Delta)_{\mathring{\Delta}\setminus A},\sd(\Delta)_{\mathring{\Delta}\setminus (A\cup \left\{F\right\})}\right)$ is a 
     good pair it holds that
     \begin{eqnarray*}
        &\empty&\widetilde{H}_{\dim \Delta}\left(\sd(\Delta)_{\mathring{\Delta}\setminus A},\sd(\Delta)_{\mathring{\Delta}
        \setminus (A\cup \left\{F\right\})};k\right)\\
        &\cong&\widetilde{H}_{\dim \Delta}\left(\sd(\Delta)_{\mathring{\Delta}\setminus A}/
        \sd(\Delta)_{\mathring{\Delta}\setminus (A\cup \left\{F\right\})};k\right)\\
        &\cong&\widetilde{H}_{\dim \Delta}\left(\mathbb{S}^{\dim \Delta};k\right) = k \\
     \end{eqnarray*}
     The same argument shows $\widetilde{H}_{\dim \Delta-1}\left(\sd(\Delta)_{\mathring{\Delta}\setminus A},
     \sd(\Delta)_{\mathring{\Delta}\setminus (A\cup\left\{F\right\})};k\right)=0$.\\
     Analogous to the case $A=\emptyset$ we deduce from the above long exact sequence
       $$
        \widetilde{H}_{\dim \Delta-1}\left(\sd(\Delta)_{\mathring{\Delta}\setminus (A\cup \left\{F\right\})};k\right)
        \widetilde{H}_{\dim \Delta-1}\left(\sd(\Delta)_{\mathring{\Delta}\setminus A};k\right)\oplus k\neq 0
       $$
     By Proposition \ref{Hochster_2} it follows that 
       $$
        \beta_{|\mathring{\Delta}\setminus (A\cup\left\{F\right\})|-\dim\Delta,|\mathring{\Delta}\setminus (A\cup\left\{F\right\})|}
        =\beta_{\sum_{l=0}^{\dim \Delta}f_l^{\Delta}-1-|A|-\dim \Delta,\sum_{l=0}^{\dim \Delta}f_l^{\Delta}-1-|A|}\neq 0. $$
    Since $A\subseteq\Delta\setminus \Delta_F$ we have $0\leq |A|\leq \sum_{l=0}^{\dim \Delta}f_l^{\Delta}-(2^{\dim \Delta+1}-1)$. 
    Therefore $\beta_{i,i+\dim \Delta}\neq 0$ for $2^{\dim \Delta+1}-2-\dim \Delta\leq i\leq \sum_{l=0}^{\dim \Delta}(f_l^{\Delta}-1)$.
  \end{proof}

  \begin{lemma} \label{lemma: Betti_3}
     Let $\Delta$ be a simplicial complex such that $\widetilde{H}_{\dim\Delta}\left(\Delta;k\right)\neq 0$ and 
     let $\beta_{ij}$ be the bigraded Betti numbers of $k[\sd(\Delta)]$.
     Then $\beta_{i,i+\dim \Delta}\neq 0$ or $\beta_{i,i+\dim \Delta +1}\neq 0$  for 
     $$2^{\dim \Delta +1}-2-\dim \Delta \leq i \leq \displaystyle{\sum_{j=0}^{\dim\Delta}}(f_j^{\Delta}-1).$$
  \end{lemma}

  \begin{proof}
     The assumption yields 
     \begin{eqnarray*}
        \widetilde{H}_{\dim \Delta}\left(\Delta;k\right)&=&\widetilde{H}_{\dim\Delta}\left(\sd(\Delta);k\right)\\
        &=&\widetilde{H}_{\sum_{l=0}^{\dim \Delta}f_l^{\Delta}-(\sum_{l=0}^{\dim\Delta}f_l^{\Delta}-\dim \Delta-1)-1}                                       
        \left(\sd(\Delta)_{\mathring{\Delta}};k\right) \neq 0.\\
     \end{eqnarray*}
     By Corollary \ref{Betti_0} it follows that $\beta_{\sum_{l=0}^{\dim\Delta}(f_l^{\Delta}-1),\sum_{l=0}^{\dim \Delta}f_l^{\Delta}}\neq 0$ 
     which proves the assertion for $i=\sum_{j=0}^{\dim \Delta}(f_j^{\Delta}-1)$.\\
     Now assume $i < \sum_{j=0}^{\dim \Delta}(f_j^{\Delta}-1)$.
     We successively remove vertices of $\Delta$ from $\sd(\Delta)$ until the
     homology in dimension $\dim \Delta$ vanishes. 
     Let $v_1, \ldots , v_r$ be vertices of $\Delta$ such that $\widetilde{H}_{\dim \Delta} (\sd(\Delta)_{\mathring{\Delta} \setminus \{\{v_1\}, \ldots, \{v_j\} \}};k)
     \neq 0$ for $1 \leq j \leq r-1$ and $\widetilde{H}_{\dim \Delta}(\sd(\Delta)_{\mathring{\Delta} \setminus \{ \{v_1\}, \ldots, \{v_r\} \}};k) =0$.
     Therefore, by Corollary \ref{Betti_0}, 
     \begin{eqnarray*}
         &\empty&\beta_{|\mathring{\Delta} \setminus \left\{\{v_1\},\ldots,\{v_j\}\right\}|-
         \dim \Delta -1,|\mathring{\Delta} \setminus \left\{\{v_1\},\ldots,\{v_j\}\right\}|}\\
         &=&\beta_{\sum_{l=0}^{\dim \Delta}f_l^{\Delta}-j-\dim \Delta-1,\sum_{l=0}^{\dim \Delta}f_l^{\Delta}-j}\\
         &\neq& 0 \text{ for } 0\leq j \leq r-1.
     \end{eqnarray*}
     Consider the complexes $\sd(\Delta)_{\mathring{\Delta}\setminus\left\{\{v_1\},\ldots,\{v_{r-1}\}\right\}}$ and 
     $\sd(\Delta)_{\mathring{\Delta}\setminus\left\{\{v_1\},\ldots,\{v_r\}\right\}}$. By construction it holds that 
     $\widetilde{H}_{\dim \Delta}\left(\sd(\Delta)_{\mathring{\Delta}\setminus\left\{\{v_1\},\ldots \{v_r\}\right\}};k\right)=0$ and\\
     $\widetilde{H}_{\dim \Delta}\left(\sd(\Delta)_{\mathring{\Delta}\setminus\left\{\{v_1\},\ldots \{v_{r-1}\}\right\}};k\right)\neq 0$.
     Successively applying the fact that for a simplicial complex $\Delta$ on ground set $\Omega$ and a vertex $v$ of $\Delta$ we get that 
     $\sd(\Delta_{\Omega \setminus \left\{v\right\}})\simeq\sd(\Delta)_{\mathring{\Delta}\setminus \left\{\{v\}\right\}}$. It follows that 
     $\sd(\Delta_{\Omega \setminus \left\{v_1,\ldots,v_{r-1}\right\}})$ $\simeq$ $\sd(\Delta)_{\mathring{\Delta}\setminus\left\{\{v_1\},\ldots,\{v_{r-1}\}\right\}}$ 
     and $\sd(\Delta_{\Omega \setminus\left\{v_1,\ldots,v_r\right\}})$ $\simeq$ $\sd(\Delta)_{\mathring{\Delta}\setminus\left\{\{v_1\},\ldots,\{v_r\}\right\}}$.
     Since the homology of a simplicial complex is invariant under barycentric subdivision this implies
     $$\widetilde{H}_{\dim \Delta}\left(\Delta_{\Omega \setminus\left\{v_1,\ldots,v_{r-1}\right\}};k\right) \neq 0 \mbox{~and~}
     \widetilde{H}_{\dim \Delta}\left(\Delta_{\Omega \setminus\left\{v_1,\ldots,v_r\right\}};k\right) = 0.$$
     Therefore, $\Delta_{\Omega \setminus \left\{v_1,\ldots,v_{r-1}\right\}}$ contains a homology cycle in dimension $\dim \Delta$. 
     We obtain the complex $\Delta_{\Omega \setminus\left\{v_1,\ldots,v_r\right\}}$ from $\Delta_{\Omega \setminus\left\{v_1,\ldots,v_{r-1}\right\}}$ 
     by removing $v_r$ and all faces containing $v_r$. Since there is a homology cycle in dimension $\dim \Delta$ the maximal dimensional faces of 
     $\Delta_{\Omega \setminus\left\{v_1,\ldots,v_{r-1}\right\}}$ cannot have a vertex in common.
     Hence, there is at least one $(\dim \Delta)$-dimensional face in $\Delta_{\Omega \setminus \left\{v_1,\ldots,v_r\right\}}$.
     Thus we have $\dim  (\Delta_{\Omega \setminus\left\{v_1,\ldots,v_r\right\}})=\dim \Delta$.\\
     Choose $F \in \Delta_{\Omega \setminus \left\{v_1,\ldots,v_r\right\}}$ with $\dim F=\dim \Delta$. 
     It follows by $\sd(\Delta_{\Omega \setminus \left\{v_1,\ldots,v_r\right\}})$ $\subseteq$ $\sd(\Delta)_{\mathring{\Delta}\setminus\{\{v_1\},\ldots,\{v_r\}\}}$ 
     that $\sd(\partial F)\subseteq \sd(\Delta)_{\mathring{\Delta}\setminus \{\{v_1\},\ldots,\{v_r\}\}}$ and $\sd(\partial F)$ is 
     homeomorphic to a $(\dim \Delta -1)$-sphere.\\
     The same arguments as in proof of Lemma \ref{lemma: Betti_2} show that
      $$\widetilde{H}_{\dim \Delta -1}\left(\sd(\Delta)_{\mathring{\Delta}\setminus(\left\{\{v_1\},\ldots,\{v_r\}\right\}\cup A\cup \left\{F\right\})};k\right)\neq 0$$
     for $A\subseteq\Delta\setminus(\Delta_F\cup\left\{\{v_1\},\ldots,\{v_r\}\right\})$. By Corollary \ref{Betti_0} 
     and $$|\mathring{\Delta}\setminus (\left\{\{v_1\},\ldots,\{v_r\}\right\}\cup A\cup \left\{F\right\})|=\sum_{l=0}^{\dim \Delta}f_l^{\Delta}-r-|A|-1$$ we deduce that 
      $$\beta_{\sum_{l=0}^{\dim \Delta}(f_l^{\Delta}-1)-r-|A|,\sum_{l=0}^{\dim\Delta}f_l^{\Delta}-r-|A|-1}\neq 0.$$
     Since $$0 \leq |A|\stackrel{A=\Delta\setminus (\Delta_F\cup\left\{v_1,\ldots,v_r\right\})}{\leq}\sum_{l=0}^{\dim\Delta} f_l^{\Delta}-r-2^{\dim \Delta +1}+1$$ 
     it follows that $\beta_{i,i+\dim\Delta}\neq 0$ for $2^{\dim \Delta +1}-\dim \Delta-2\leq i\leq \sum_{l=0}^{\dim\Delta}(f_l^{\Delta}-1)-r$ what finally completes the proof.
  \end{proof}

  The following lemma is a simple consequence of the characterization \cite[Theorem 1]{BK} 
  of pairs of
  the vector $(\dim \widetilde{H}_i(\Delta;k))_{0 \leq i \leq \dim \Delta}$ encoding the
  Betti numbers of $\Delta$ and the $f$-vector $(f_i^\Delta)_{-1 \leq i \leq \dim \Delta}$ 
  of $\Delta$. We leave the verification to the reader.

  \begin{lemma} \label{lemma: f-vector}
     Let $\Delta$ be a $d$-dimensional simplicial complex such that $\widetilde{H}_d(\Delta;k)=0$.
     Then $f_{d-1}^{\Delta}\geq f_d^{\Delta}+d$.
  \end{lemma}

  The following lemmas include simple but crucial inequalities that will be 
  used in the derivation of the main theorem. Their proofs are straightforward
  and left to the reader. 

  \begin{lemma} \label{lemma no2}
     For $d\geq 1$ 
     $$\frac{\prod_{l=2}^{d+1}(2^{d+1}-l)}{(d+1)!\cdot \prod_{m=2}^d(2^{m+1}-3)}\geq 
     \begin{cases}
        1\text{ if }1\leq d\leq 3\\
        2 \text{ if } d\geq 4\\
     \end{cases}$$
  \end{lemma}

  \begin{lemma}\label{lemma no3}
     For $n\geq 11$ 
        $$(n+1)!\leq 2^{\frac{n^2}{2}-\frac{5}{2}n}.$$
  \end{lemma}

  \begin{lemma}\label{lemma no4}
     For $n\in \mathbb{N}$ and $k\geq 2$ 
     $$\frac{\prod_{l=0}^{n-1}(2^{n+1}+2k-4+l)}{(n+1)!\cdot k\cdot 
     \prod_{m=2}^{n}(2^{m+1}-3)}\geq 1.$$
  \end{lemma}

  \begin{lemma} \label{lemma no5}
     For $d\geq 4$ it holds that 
     $$d \cdot\prod_{l=0}^{d-2}(2^{d+2}-d-6-l)\geq (d+1)!\cdot \prod_{l=2}^{d}(2^{l+1}-3).$$
  \end{lemma}

\section{Proof of the Theorem \ref{MainResult}} \label{ProofOfMainResult}

   Before we proceed to the proof of the main theorem, we consider $\Delta$ with small dimension.
   If $\dim \Delta = 0$ then $I_\Delta$ is generated by all squarefree monomials of degree $2$.
   It is well known that the resolution of this ideal is linear. Hence the Multiplicity Conjecture
   holds (see e.g. \cite{HS}).
   For the cases $\dim \Delta =1,2$ the conjecture was settled in \cite[Theorem 4.3]{NS} except for the
   equality statement.
   For dimensions $3$ and $4$ the result follows from
   \cite{NS} in case the complex is Gorenstein.

   \begin{proof}[Proof of Theorem \ref{MainResult}]
      {\sf Upper Bound:} By the argumentation above we may assume that $\dim(\Delta) = \dim(\sd(\Delta)) \geq 1$.
      We set $F^\Delta := \sum_{l=0}^{\dim \Delta}(f_l^{\Delta}-1)$.
      By Proposition \ref{height, multiplicity} we have to show that 
      $$(\dim \Delta +1)!\cdot f_{\dim \Delta}^{\Delta}\leq\frac{1}
      {\left(F^\Delta\right)!}\cdot 
      \prod_{i=1}^{F^\Delta}M_i.$$

      First, we consider the case $\widetilde{H}_{\dim \Delta}\left(\Delta;k\right)=0$. 
      From Lemmas \ref{lemma: Betti_1} and \ref{lemma: Betti_2} we deduce that 
      $M_i\geq m+i$ for  $2^{m+1}-2-m\leq i< 2^{m+2}-2-(m+1)$ and $1\leq m < \dim \Delta$ 
      and $M_i\geq i+\dim \Delta$ for $2^{\dim \Delta +1}-2-
      \dim \Delta\leq i\leq F^\Delta$. Therefore:
      \begin{eqnarray*}
         &\empty&\prod_{i=1}^{F^\Delta}M_i\\
         &=&\prod_{m=1}^{\dim \Delta -1}\left(\prod_{i=2^{m+1}-2-m}^{2^{m+2}-2-(m+1)-1}M_i\right)
         \cdot \prod_{i=2^{\dim \Delta+1}-2-\dim \Delta}^{F^\Delta}M_i\\
         &\geq&\prod_{m=1}^{\dim \Delta-1}\left(\prod_{i=2^{m+1}-2-m}^{2^{m+2}-m-4}(m+i)\right)
         \cdot \prod_{i=2^{\dim \Delta+1}-2-\dim \Delta}^{F^\Delta}(i+\dim \Delta)\\
         &=&\left(\prod_{m=1}^{\dim \Delta -1}\frac{\left(2^{m+2}-m-4+m\right)!}
         {\left(2^{m+1}-2-m+m-1\right)!}\right)\\
	 &\empty&\cdot \frac{\left(F^\Delta+\dim \Delta\right)!}
         {\left(2^{\dim \Delta +1}-2-\dim \Delta+\dim \Delta-1\right)!}\\
         &=& \left(\prod_{m=1}^{\dim \Delta -1}\frac{\left(2^{m+2}-4\right)!}
         {\left(2^{m+1}-3\right)!}\right)\cdot \frac{(F^\Delta+\dim\Delta)!}{\left(2^{\dim \Delta +1}-3\right)!}\\
         &=& \frac{\prod_{m=2}^{\dim \Delta}\left(2^{m+1}-4\right)!}
         {\prod_{m=1}^{\dim \Delta}\left(2^{m+1}-3\right)!}\cdot 
         (F^\Delta+\dim\Delta)!\\
         &=& \frac{\prod_{m=2}^{\dim \Delta}\left(2^{m+1}-4\right)!}
         {\prod_{m=2}^{\dim \Delta}\left(2^{m+1}-3\right)!}\cdot 
         \left(F^\Delta+\dim \Delta\right)!\\
         &=& \prod_{m=2}^{\dim \Delta}\frac{1}{2^{m+1}-3}\cdot 
         \left(F^\Delta+\dim \Delta\right)!.
      \end{eqnarray*}
      It follows
      \begin{eqnarray*}
         &\empty&\frac{1}{\left(F^\Delta\right)
         !\cdot (\dim \Delta +1)!\cdot f_{\dim \Delta}^{\Delta}}\cdot 
         \prod_{i=1}^{F^\Delta}M_i\\
         &\geq&\frac{1}{\left(F^\Delta\right)
         !\cdot (\dim \Delta+1)!\cdot f_{\dim \Delta}^{\Delta}}\cdot 
         \prod_{m=2}^{\dim\Delta}\frac{1}{2^{m+1}-3}\cdot 
         \left(F^\Delta+\dim \Delta\right)!\\
         &=&\frac{\prod_{m=1}^{\dim \Delta}\left(F^\Delta+\dim \Delta+1-m\right)}{(\dim \Delta +1)!\cdot f_{\dim \Delta}^{\Delta}\cdot 
         \prod_{m=2}^{\dim \Delta}(2^{m+1}-3)}.
       \end{eqnarray*}
       Along with Lemma \ref{lemma: f-vector} this yields
       \begin{eqnarray*}
          &\empty&\frac{1}{\left(F^\Delta\right)
          !\cdot (\dim \Delta +1)!\cdot f_{\dim \Delta}^{\Delta}}\cdot 
          \prod_{i=1}^{F^\Delta}M_i\\
          &\geq&\frac{\prod_{m=1}^{\dim \Delta}\left(\sum_{l=0}^{\dim \Delta-2}
          f_l^{\Delta}+2f_{\dim \Delta}^{\Delta}+\dim \Delta-m\right)}{(\dim \Delta +1)!
          \cdot f_{\dim \Delta}^{\Delta}\cdot \prod_{m=2}^{\dim \Delta}(2^{m+1}-3)}\\
          &=&\frac{\prod_{m=0}^{\dim \Delta-1}\left(\sum_{l=0}^{\dim \Delta-2}
          f_l^{\Delta}+2f_{\dim \Delta}^{\Delta}+m\right)}{(\dim \Delta +1)!\cdot 
          f_{\dim \Delta}^{\Delta}\cdot \prod_{m=2}^{\dim \Delta}(2^{m+1}-3)}.
       \end{eqnarray*}
       Assume that $f_{\dim \Delta}^{\Delta}=1$. Then $f_i^{\Delta}\geq 
       \binom{\dim \Delta +1}{i+1}$, $0\leq i\leq \dim \Delta$. This implies
       \begin{eqnarray*}
          &\empty&\frac{1}{\left(F^\Delta\right)!\cdot 
          (\dim \Delta +1)!\cdot f_{\dim \Delta}^{\Delta}}\cdot 
          \prod_{i=1}^{F^\Delta}M_i\\
          &\geq&\frac{\prod_{m=0}^{\dim \Delta -1}\left(\sum_{l=0}^{\dim \Delta -2}
          \binom{\dim \Delta+1}{l+1}+2+m\right)}{(\dim \Delta+1)!\cdot 
          \prod_{m=2}^{\dim \Delta}(2^{m+1}-3)}\\
          &=&\frac{\prod_{m=0}^{\dim \Delta-1}\left(2^{\dim \Delta+1}-
          (1+ \binom{\dim\Delta+1}{\dim \Delta}+1)+2+m\right)}{(\dim \Delta+1)!\cdot 
          \prod_{m=2}^{\dim \Delta}(2^{m+1}-3)}\\
          &=&\frac{\prod_{m=0}^{\dim \Delta-1}\left(2^{\dim \Delta+1}-
          \dim\Delta-1+m\right)}{(\dim \Delta+1)!\cdot \prod_{m=2}^{\dim \Delta}(2^{m+1}-3)}\\
          &=&\frac{\prod_{m=2}^{\dim \Delta+1}\left(2^{\dim \Delta+1}-m\right)}
          {(\dim \Delta+1)!\cdot \prod_{m=2}^{\dim \Delta}(2^{m+1}-3)}.
       \end{eqnarray*}
       Since by Lemma \ref{lemma no2} the latter expression is greater or equal than $1$ 
       this shows the claim in case $\widetilde{H}_{\dim \Delta}\left(\Delta;k\right)=0$ and 
       $f_{\dim \Delta}^{\Delta}=1$.\\
       Let $f_{\dim \Delta}^{\Delta}>1$. Clearly, in this case $f_i^{\Delta}\geq \binom{\dim \Delta+1}{i+1}+1$. 
       Therefore
       \begin{eqnarray*}
          &\empty&\frac{1}{\left(F^\Delta\right)!\cdot 
          (\dim \Delta +1)!\cdot f_{\dim \Delta}^{\Delta}}\cdot 
          \prod_{i=1}^{F^\Delta}M_i\\
          &\geq&\frac{\prod_{m=0}^{\dim \Delta -1}\left(\sum_{l=0}^{\dim \Delta -2}
          (\binom{\dim \Delta+1}{l+1}
          +1)+2f_{\dim \Delta}^{\Delta}+m\right)}{(\dim \Delta+1)!\cdot 
          f_{\dim \Delta}^{\Delta}\cdot \prod_{m=2}^{\dim \Delta}(2^{m+1}-3)}\\
          &=&\frac{\prod_{m=0}^{\dim \Delta -1}\left(\sum_{l=0}^{\dim \Delta -2}
          \binom{\dim \Delta+1}{l+1}+\dim \Delta-1+2f_{\dim \Delta}^{\Delta}+m\right)}
          {(\dim \Delta+1)!\cdot f_{\dim \Delta}^{\Delta}\cdot 
          \prod_{m=2}^{\dim \Delta}(2^{m+1}-3)}\\
          &=&\frac{\prod_{m=0}^{\dim \Delta -1}\left(2^{\dim \Delta+1}-
          (1+(\dim \Delta+1)+1)+\dim \Delta-1+2f_{\dim \Delta}^{\Delta}+m\right)}
          {(\dim \Delta+1)!\cdot f_{\dim \Delta}^{\Delta}\cdot 
          \prod_{m=2}^{\dim \Delta}(2^{m+1}-3)}\\
          &=&\frac{\prod_{m=0}^{\dim \Delta-1}(2^{\dim \Delta+1}-4+
          2f_{\dim \Delta}^{\Delta}+m)}{(\dim \Delta+1)!\cdot 
          f_{\dim \Delta}^{\Delta}\cdot \prod_{m=2}^{\dim \Delta}(2^{m+1}-3)}\geq 1.
       \end{eqnarray*}
       The last inequality holds by Lemma \ref{lemma no4}. This proves the 
       Multiplicity Conjecture if $\widetilde{H}_{\dim \Delta}\left(\Delta;k\right)=0$ and 
       $f_{\dim \Delta}^{\Delta}>1$.\\
       Let $\widetilde{H}_{\dim \Delta}\left(\Delta;k\right)\neq 0$.
       If $\Delta$ has dimension $1$ or $2$ the claim follows from \cite[Theorem 4.3]{NS}.
       Let $\dim \Delta \geq 3$. By Lemma \ref{lemma: Betti_1} and \ref{lemma: Betti_3} it 
       holds that $M_i\geq m+i$ for $1\leq m < \dim \Delta$ and $2^{m+1}-2-m\leq i< 2^{m+2}-2-(m+1)$ 
       and $M_i\geq i+\dim \Delta$ for $2^{\dim \Delta +1}-2-\dim \Delta\leq i\leq 
       F^\Delta-1$ and $M_{F^\Delta}\geq F^\Delta+
       \dim \Delta+1=\sum_{l=0}^{\dim \Delta}f_l^{\Delta}$. Therefore, the same 
       calculation as in the first part of the proof yields
       $$\prod_{i=1}^{F^\Delta}M_i
       \geq\prod_{m=2}^{\dim \Delta}\frac{1}{2^{m+1}-3}\cdot 
       \frac{\left(F^\Delta+\dim \Delta+1\right)!}
       {F^\Delta+\dim \Delta}.$$
       Thus it suffices to show that 
       $$\frac{\left(F^\Delta+\dim\Delta+1\right)!}{(\dim \Delta+1)! 
       f_{\dim \Delta}^{\Delta}\left(F^\Delta\right)
       !\left(F^\Delta+\dim \Delta\right)\prod_{m=2}^{\dim \Delta}
       (2^{m+1}-3)}\geq 1.$$
       We have that 
       \begin{eqnarray*}
          &\empty&\frac{\left(F^\Delta+\dim \Delta+1\right)!}{(\dim \Delta+1)! 
          f_{\dim \Delta}^{\Delta}\left(F^\Delta\right)
          !\left(F^\Delta+\dim \Delta\right)\prod_{m=2}^{\dim \Delta}
          (2^{m+1}-3)}\\
          &=&\frac{\prod_{m=2}^{\dim \Delta}\left(F^\Delta+\dim \Delta+1-m\right)\cdot (F^\Delta+\dim \Delta+1)}
          {(\dim \Delta+1)!f_{\dim \Delta}^{\Delta}\prod_{m=2}^{\dim \Delta}(2^{m+1}-3)}.
       \end{eqnarray*}
       Since $\widetilde{H}_{\dim \Delta}\left(\Delta;k\right)\neq 0$ it holds that 
       $f_i^{\Delta}\geq f_i^{\partial(\Delta_{\dim \Delta+1})}=
       \binom{\dim \Delta+2}{i+1}$, where $\partial(\Delta_{\dim \Delta+1})$ 
       denotes the boundary of the $(\dim \Delta+1)$-simplex. It follows that 
       \begin{eqnarray}
          \sum_{l=0}^{\dim \Delta-1}f_l^{\Delta}&\geq& 
          \sum_{l=0}^{\dim \Delta-1}\binom{\dim \Delta+2}{l+1}\nonumber\\
          &=&2^{\dim \Delta+2}-(1+\binom{\dim \Delta+2}{\dim \Delta+1}+1)\nonumber\\
          &=&2^{\dim \Delta+2}-\dim \Delta -4.\label{eqn}
       \end{eqnarray}
       We conclude that
       \begin{eqnarray*}
          &\empty&\frac{\left(F^\Delta+\dim \Delta+1\right)!}
          {(\dim \Delta+1)! f_{\dim \Delta}^{\Delta}
          \left(F^\Delta\right)!
          \left(F^\Delta+\dim \Delta\right)
          \prod_{m=2}^{\dim \Delta}(2^{m+1}-3)}\\
          &\geq&\prod_{m=2}^{\dim \Delta}\left(2^{\dim \Delta+2}-
          \dim \Delta-4+f_{\dim \Delta}^{\Delta}-m\right)\\
	 &\empty&\cdot\frac{(2^{\dim \Delta+2}-\dim \Delta-4+f_{\dim \Delta}^{\Delta})}{(\dim \Delta+1)!f_{\dim \Delta}^{\Delta}\prod_{m=2}^{\dim \Delta}(2^{m+1}-3)}\\
          &=&\prod_{m=6+\dim \Delta}^{4+2\dim \Delta}
          \left(2^{\dim \Delta+2}-m+f_{\dim \Delta}^{\Delta}\right) \\
	&\empty&\cdot\frac{\left(2^{\dim \Delta+2}-\dim \Delta-4+
          f_{\dim \Delta}^{\Delta}\right)}{(\dim \Delta+1)!
          f_{\dim \Delta}^{\Delta}\prod_{m=2}^{\dim \Delta}(2^{m+1}-3)}\\
          &\geq& \frac{\dim \Delta\cdot f_{\dim \Delta}^{\Delta}\cdot
          \prod_{m=6+\dim \Delta}^{4+2\dim \Delta}(2^{\dim \Delta+2}-m)} 
          {(\dim \Delta+1)!f_{\dim \Delta}^{\Delta}\prod_{m=2}^{\dim \Delta}(2^{m+1}-3)}\\
          &=&\frac{\dim \Delta\cdot\prod_{m=6+\dim \Delta}^{4+2\dim \Delta}
          (2^{\dim \Delta+2}-m)}{(\dim \Delta+1)!\prod_{m=2}^{\dim \Delta}(2^{m+1}-3)}\\
          &=&\frac{\dim \Delta\cdot\prod_{m=0}^{\dim \Delta-2}(2^{\dim \Delta+2}-
          \dim\Delta-6-m)}{(\dim \Delta+1)!\prod_{m=2}^{\dim \Delta}(2^{m+1}-3)}\\
          &\geq& 1,
       \end{eqnarray*}
       where the last inequality holds by Lemma \ref{lemma no5} for $\dim \Delta\geq4$.\\
       It remains to show the assertion for $\dim \Delta=3$.\\
       By Equation (\ref{eqn}) we have that 
       \begin{eqnarray*}
          &\empty&\frac{(\sum_{l=0}^{3}f_l^{\Delta}-3)(\sum_{l=0}^3f_l^{\Delta}-2)
          (\sum_{l=0}^{3}f_l^{\Delta})}{4!\cdot f_3^{\Delta}\cdot \prod_{m=2}^{3}(2^{m+1}-3)}\\
          &\geq&\frac{(2^5-10+f_3^{\Delta})(2^5-9+f_3^{\Delta})(2^5-7+f_3^{\Delta})}
          {24\cdot f_3^{\Delta}\cdot 5 \cdot 13}\\
          &=&\frac{(22+f_3^{\Delta})(23+f_3^{\Delta})(25+f_3^{\Delta})}{1560f_3^{\Delta}}\\
          &=&\frac{(f_3^{\Delta})^3+70(f_3^{\Delta})^2+1631f_3^{\Delta}+12650}{1560f_3^{\Delta}}\\
          &\geq&\frac{1631 f_3^{\Delta}}{1560f_3^{\Delta}}\geq 1
       \end{eqnarray*}
       This finally concludes the proof of the upper bound in the Multiplicity Conjecture.\\

      \noindent{\sf Cohen-Macaulay Case and Lower Bound:}
      By Proposition \ref{height, multiplicity} we have to show that 
      $$\frac{1}{\left(F^\Delta\right)!}\cdot 
      \prod_{i=1}^{F^\Delta}m_i\leq(\dim \Delta +1)!\cdot f_{\dim \Delta}^{\Delta}.$$
     \noindent $\triangleright$ {\sf Case 1:} First, we consider simplicial complexes $\Delta$ for which $\dim \Delta \geq 7$ and there exists a 
     face $F \in \Delta$ of dimension $\dim\Delta-1$ that is contained in a unique facet $G$. 
     Then, the restricted complex $\sd(\Delta)_{\mathring{\Delta}_G\setminus(\{G\}\cup\mathring{\partial F})}$
     consists of at least two connected components. Therefore, $\sd(\Delta)_W$ is disconnected if $|W|\geq 2$,
     $F\in W$ and $W\subseteq \mathring{\Delta}\setminus(\{G\}\cup\mathring{\partial F})$.
     This implies $\widetilde{H}_{0}(\sd(\Delta)_W;k)\neq 0$ for $|W|\geq 2$,
     $F\in W$ and $W\subseteq \mathring{\Delta}\setminus(\{G\}\cup\mathring{\partial F})$.
     From Corollary \ref{Betti_0} we deduce $\beta_{|W|-1,|W|}\neq 0$ for $2\leq |W|\leq \sum_{j=0}^{\dim\Delta}f_j^{\Delta}-(2^{\dim\Delta}-2)-1$,
     i.e. $\beta_{i,i+1}\neq 0$ for $1\leq |W|\leq \sum_{j=0}^{\dim\Delta}f_j^{\Delta}-2^{\dim\Delta}$. Thus $m_i\leq i+1$ for
     $1\leq |W|\leq \sum_{j=0}^{\dim\Delta}f_j^{\Delta}-2^{\dim\Delta}$. From Section \ref{Regularity} we know that $\reg(k[\sd(\Delta)])\leq \dim\Delta+1$ . 
     Hence, $m_i\leq i+\dim \Delta+1$ for $ \sum_{j=0}^{\dim\Delta}f_j^{\Delta}-2^{\dim\Delta}+1\leq i\leq F^{\Delta}$.
     This implies 
     \begin{eqnarray*}
     &\empty&\prod_{i=1}^{F^{\Delta}}m_i\cdot\frac{1}{(F^{\Delta})!}\\
     &=&\prod_{i=1}^{\sum_{j=0}^{\dim\Delta}f_j^{\Delta}-2^{\dim\Delta}}m_i\cdot\prod_{i=\sum_{j=0}^{\dim\Delta}f_j^{\Delta}-2^{\dim\Delta}+1}^{F^{\Delta}}m_i\cdot\frac{1}{(F^{\Delta})!}\\
     &\leq&\frac{1}{(F^{\Delta})!}\cdot\prod_{i=1}^{\sum_{j=0}^{\dim\Delta}f_j^{\Delta}-2^{\dim\Delta}}(i+1)\cdot\prod_{i=\sum_{j=0}^{\dim\Delta}f_j^{\Delta}-2^{\dim\Delta}+1}^{F^{\Delta}}(i+\dim\Delta+1)\\
     &=&\frac{1}{(F^{\Delta})!}\cdot(\sum_{j=0}^{\dim\Delta}f_j^{\Delta}-2^{\dim\Delta}+1)!\cdot\frac{(F^{\Delta}+\dim\Delta+1)!}{(\sum_{j=0}^{\dim\Delta}f_j^{\Delta}-2^{\dim\Delta}+\dim\Delta+1)!}\\
     &=& \frac{\prod_{i=1}^{\dim\Delta+1}(F^{\Delta}+i)}{\prod_{i=1}^{\dim\Delta}(\sum_{j=0}^{\dim\Delta}f_j^{\Delta}-2^{\dim\Delta}+1+i)}.
     \end{eqnarray*}
     The claim now follows if we show that
    \begin{eqnarray}
    &\empty&\prod_{i=1}^{\dim\Delta+1}(F^{\Delta}+i)\label{Beh_1}\\
    &\leq& (\dim \Delta+1)!\cdot f_{\dim\Delta}^{\Delta}\cdot\prod_{i=1}^{\dim\Delta}(\sum_{j=0}^{\dim\Delta}f_j^{\Delta}-2^{\dim\Delta}+1+i).\nonumber
    \end{eqnarray}
     For $d\geq8$ it holds that $d!\geq 2^{2d-1}$. From 
     \begin{eqnarray*}
     2^{\dim\Delta+1}\cdot f_{\dim\Delta}^{\Delta}&\geq& \sum_{j=0}^{\dim\Delta}f_j^{\Delta}\\
     &=&F^{\Delta}+\dim\Delta+1
     \end{eqnarray*}
     we conclude that
     \begin{eqnarray*}
     (\dim\Delta+1)!\cdot f_{\dim\Delta}^{\Delta}&\geq& 2^{2\cdot\dim\Delta+1}\cdot f_{\dim\Delta}^{\Delta}\\
     &=&2^{\dim\Delta}\cdot 2^{\dim\Delta+1}\cdot f_{\dim\Delta}^{\Delta}\\
     &\geq& 2^{\dim\Delta}\cdot(F^{\Delta}+\dim\Delta+1)
     \end{eqnarray*}
     for $\dim\Delta\geq 7$.\\
     Therefore, it suffices to show that
    $$\prod_{i=1}^{\dim\Delta+1}(F^{\Delta}+i)\leq 2^{\dim\Delta}\cdot (F^{\Delta}+\dim\Delta+1)\cdot\prod_{i=1}^{\dim\Delta}(\sum_{j=0}^{\dim\Delta}f_j^{\Delta}-2^{\dim\Delta}+1+i),$$
    i.e.
    $$\prod_{i=1}^{\dim\Delta}(F^{\Delta}+i)\leq\prod_{i=1}^{\dim\Delta}2\cdot(\sum_{j=0}^{\dim\Delta}f_j^{\Delta}-2^{\dim\Delta}+1+i).$$
    By $\sum_{j=-1}^{\dim\Delta}f_j^{\Delta}\geq 2^{\dim\Delta+1}$ it follows that 
    \begin{eqnarray*}
    &\empty&F^{\Delta}+2\dim\Delta+2+i\\
    &=&\sum_{j=-1}^{\dim\Delta}f_j^{\Delta}+\dim\Delta+1+i\\
    &\geq&2^{\dim\Delta+1} 
     \end{eqnarray*}
     for $1\leq i\leq \dim\Delta$. This implies 
     \begin{eqnarray*}
     &\empty&F^{\Delta}+i\\
     &\leq& 2\cdot(F^{\Delta}+\dim\Delta+1-2^{\dim\Delta}+1+i)\\
     &=& 2\cdot(\sum_{j=0}^{\dim\Delta}f_j^{\Delta}-2^{\dim\Delta}+1+i)
     \end{eqnarray*}
     for $1\leq i\leq \dim\Delta$. This concludes the proof in Case 1. \\
     \noindent $\triangleright$ {\sf Case 2:} 
     Now, we consider simplicial complexes $\Delta$ such that $\dim \Delta \geq 6$ and every $(\dim\Delta-1)$-dimensional face lies in at 
     least two facets.
     Let $F \in \Delta$ be a face of dimension $\dim\Delta$. The restriction $\sd(\Delta)_{\mathring{\Delta}\setminus\mathring{\partial F}}$ 
     consists of two connected components. In particular, every restriction of the form $\sd(\Delta)_{W}$ where 
     $W\subseteq\mathring{\Delta}\setminus\mathring{\partial F}$, $F\in W$ and $|W|\geq 2$, is disconnected, 
     i.e. $\widetilde{H}_0(\sd(\Delta)_W;k)\neq 0$. From Corollary \ref{Betti_0} we deduce that $\beta_{|W|-1,|W|}\neq 0$. Therefore, 
     $\beta_{i,i+1}\neq 0$ for $2\leq i+1\leq \sum_{j=0}^{\dim\Delta}f_j^{\Delta}-(2^{\dim\Delta+1}-2)$. This implies
     $m_i\leq i+1$ for $1\leq i\leq \sum_{j=0}^{\dim\Delta}f_j^{\Delta}-2^{\dim\Delta+1}+1$. As in Case 1, Section \ref{Regularity} implies that 
     $m_i\leq i+\dim\Delta+1$ for $\sum_{j=0}^{\dim\Delta}f_j^{\Delta}-2^{\dim\Delta+1}+2\leq i\leq F^{\Delta}$. Therefore, 

     \begin{eqnarray*}
     &\empty&\frac{1}{(F^{\Delta})!}\cdot \prod_{i=1}^{F^{\Delta}}m_i\\
     &=&\frac{1}{(F^{\Delta})!}\cdot\prod_{i=1}^{\sum_{j=0}^{\dim\Delta}f_j^{\Delta}-2^{\dim\Delta+1}+1}m_i\cdot\prod_{i=\sum_{j=0}^{\dim\Delta}f_j^{\Delta}-2^{\dim\Delta+1}+2}^{F^{\Delta}}m_i\\
     &\leq&\frac{1}{(F^{\Delta})!}\cdot\prod_{i=1}^{\sum_{j=0}^{\dim\Delta}f_j^{\Delta}-2^{\dim\Delta+1}+1}(i+1)\cdot\prod_{i=\sum_{j=0}^{\dim\Delta}f_j^{\Delta}-2^{\dim\Delta+1}+2}^{F^{\Delta}}(i+\dim\Delta+1)\\
     &=&\frac{(F^{\Delta}+\dim\Delta+1-2^{\dim\Delta+1}+2)!}{(F^{\Delta})!}\cdot\frac{(F^{\Delta}+\dim\Delta+1)!}{(F^{\Delta}+2\dim\Delta+2-2^{\dim\Delta+1}+1)!}\\
     &=&\frac{\prod_{i=1}^{\dim\Delta+1}(F^{\Delta}+i)}{\prod_{i=1}^{\dim\Delta}(F^{\Delta}+\dim\Delta+1-2^{\dim\Delta+1}+2+i)}.
     \end{eqnarray*}
     Thus, it suffices to show that
     $$\prod_{i=1}^{\dim\Delta+1}(F^{\Delta}+i)\leq(\dim\Delta+1)!\cdot f_{\dim\Delta}^{\Delta}\cdot\prod_{i=1}^{\dim\Delta}(F^{\Delta}+\dim\Delta+3-2^{\dim\Delta+1}+i).$$
     Since $\Delta$ is pure, every $i$-dimensional face is contained in a $(\dim\Delta-1)$-dimensional face and then by assumption in at least two 
     facets of $\Delta$. This implies
     $f_i^{\Delta}\leq \binom{\dim\Delta+1}{ i+1}\cdot f_{\dim\Delta}^{\Delta}\cdot\frac{1}{2}$ for $0\leq i\leq \dim\Delta-1$. Hence
     \begin{eqnarray}
     \sum_{i=0}^{\dim\Delta}f_i^{\Delta}&\leq&\sum_{i=0}^{\dim\Delta-1}\binom{\dim\Delta+1}{i+1}\cdot f_{\dim\Delta}^{\Delta}\cdot\frac{1}{2}+f_{\dim\Delta}^{\Delta}\nonumber \\
     &=&(2^{\dim\Delta+1}-2)\cdot f_{\dim\Delta}^{\Delta}\cdot\frac{1}{2}+f_{\dim\Delta}^{\Delta}\nonumber\\
     &=&2^{\dim\Delta}\cdot f_{\dim\Delta}^{\Delta},\label{facenumbers}
     \end{eqnarray}

     i.e. $2^{\dim\Delta}\cdot f_{\dim\Delta}^{\Delta}\geq F^{\Delta}+\dim\Delta+1$.\\
     Using $d!\geq 2^{2d-2}$ for $d\geq 7$ we obtain that 
     $(\dim\Delta+1)!\cdot f_{\dim\Delta}^{\Delta}\geq 2^{\dim\Delta}\cdot(F^{\Delta}+\dim\Delta+1)$ for $\dim\Delta\geq 6$. Therefore, it suffices to show that 
     $$\prod_{i=1}^{\dim\Delta+1}(F^{\Delta}+i)\leq 2^{\dim\Delta}\cdot(F^{\Delta}+\dim\Delta+1)\cdot\prod_{i=1}^{\dim\Delta}(F^{\Delta}+\dim\Delta+3-2^{\dim\Delta+1}+i),$$
     i.e.
     $$\prod_{i=1}^{\dim\Delta}(F^{\Delta}+i)\leq\prod_{i=1}^{\dim\Delta}2\cdot(F^{\Delta}+\dim\Delta+3-2^{\dim\Delta+1}+i)$$
     for $\dim\Delta\geq 6$.\\
     For a simplicial complex with the property that every face of dimension $\dim\Delta-1$ is contained in at least two facets each entry of its
     $f$-vector is bounded from below by the same entry in the $f$-vector of the boundary of the $(\dim\Delta+1)$-simplex. Hence, 
     \begin{eqnarray}
     \sum_{j=-1}^{\dim\Delta}f_j^{\Delta}\geq 2^{\dim\Delta+2}-1.\label{face_number}
     \end{eqnarray}
      This implies
     \begin{eqnarray*}
     &\empty&F^{\Delta}+2\cdot\dim\Delta+6+i\\
     &=&\sum_{j=-1}^{\dim\Delta}f_j^{\Delta}+\dim\Delta+4+i\geq 2^{\dim\Delta+2}
     \end{eqnarray*}
     for $1\leq i\leq\dim\Delta$. We deduce that 
     $$F^{\Delta}+i\leq 2\cdot (F^{\Delta}+\dim\Delta+3-2^{\dim\Delta+1}+i)$$
      for $1\leq i\leq \dim\Delta$.  This proves the lower bound conjecture in Case 2.

     $\triangleright$ {\sc Small Dimensions:} 
     We are now the lower bound of the Multiplicity Conjecture for the dimensions not covered in Case 1 and Case 2. \\
     Since $\Delta$ is a Cohen-Macaulay complex there exists an ordering $F_1,\ldots, F_{f_{\dim\Delta}^{\Delta}}$ of the facets of $\Delta$ such that
     $\dim\overline{(F_1\cup\ldots\cup F_i)}\cap F_{i+1}=\dim\Delta-1$. This implies
     \begin{eqnarray}
     &\empty&\sum_{j=0}^{\dim\Delta}f_j^{\Delta}\nonumber\\ 
     &\leq& (2^{\dim\Delta+1}-1)+(f_{\dim\Delta}^{\Delta}-1)\cdot (2^{\dim\Delta+1}-1-(2^{\dim\Delta}-1)) \nonumber\\
     &=& 2^{\dim\Delta}-1+2^{\dim\Delta}\cdot f_{\dim\Delta}^{\Delta}.\label{ordering}
     \end{eqnarray}
     By the arguments preceding the proof we may assume $\dim \Delta \geq 3$.

     ${\mathbf \dim\Delta=3}.$ In the situation of Case 1 we have to show that 
     $$(F^{\Delta}+1)\cdots  (F^{\Delta}+4)\leq 4!\cdot f_3^{\Delta}\cdot (F^{\Delta}+4-2^3+1+1)\cdot (F^{\Delta}-1)\cdot F^{\Delta}.$$
     Inequality (\ref{ordering}) yields $F^{\Delta}+4\leq 8\cdot f_3^{\Delta}+7$. If $\Delta$ is the $3$-simplex, then $\sd(\Delta)$ is
     Gorenstein an the result follows by the \cite{NS}.
% some calculations show that
%     $m_i\leq i+1$ for $1\leq i\leq 7$, as well as $m_i\leq i+2$ for $8\leq i\leq 10$ and $m_{11}=14$. 
%     Since $e(k[\sd(\Delta)])=96$ and $F^{\Delta}=11$ we get 
%     $\frac{1}{11!}\cdot \prod_{i=1}^{11}m_i\leq \frac{56}{3}\leq 96 =e(k[\sd(\Delta)])$.\\
     If $\Delta$ is not the $3$-simplex it follows from the above considerations that it suffices to show that
     $$(F^{\Delta}+1)\cdot (F^{\Delta}+2)\cdot (F^{\Delta}+3)\leq 2\cdot (F^{\Delta}-2)\cdot (F^{\Delta}-1)\cdot F^{\Delta}$$
      which is equivalent to $(F^{\Delta})^3-12\cdot (F^{\Delta})^2-7\cdot F^{\Delta}-6\geq 0$. Since $\Delta$ is not the $3$-simplex we can
      deduce $F^{\Delta}\geq 2^4-1+(2^4-1-(2^3-1))-4=19$. Since 
     \begin{eqnarray*}
     6+7\cdot F^{\Delta}+12\cdot (F^{\Delta})^2&\leq& 8\cdot F^{\Delta}+12\cdot (F^{\Delta})^2\\
     &\leq& (F^{\Delta})^2+12\cdot (F^{\Delta})^2\\
     &=& 13\cdot (F^{\Delta})^2\leq (F^{\Delta})^3
      \end{eqnarray*}
      for $F^{\Delta}\geq 19$ the claim follows.\\
      Now we turn to the situation of Case 2. If $\Delta$ is the boundary of the $4$-simplex then the complex and its barycentric subdivision
      are Gorenstein. Hence the result follows from \cite{NS}. Assume $\Delta$ is not the boundary of the $4$-simplex.  Then $\Delta$ must 
      have at least one additional $3$-simplex (i.e, $f_3^\Delta \geq 5$). From the Kruskal-Katona theorem we infer that
      $f_2^\Delta,f_1^\Delta \geq 12$ and $f_0^\Delta \geq 6$. This implies $F^\Delta \geq 32$. We have to show that
      $$(F^{\Delta}+1)\cdot \ldots\cdot (F^{\Delta}+4)\leq 4!\cdot f_3^{\Delta}\cdot (F^{\Delta}+4-2^4+2+1)\cdot (F^{\Delta}-8)\cdot (F^{\Delta}-7).$$
      By Inequality (\ref{facenumbers}) it holds that $F^{\Delta}+4\leq 8\cdot f_3^{\Delta}$. Hence, it suffices to show that
     $(F^{\Delta}+1)\cdot (F^{\Delta}+2)\cdot (F^{\Delta}+3)\leq 3\cdot (F^{\Delta}-9)\cdot (F^{\Delta}-8)\cdot (F^{\Delta}-7)$. 
% , i.e.  $(F^{\Delta})^3-39\cdot (F^{\Delta})^2+281\cdot F^{\Delta}-759\geq 0$.  
     The inequality is satisfied for $F^{\Delta}\geq 31$ and we are done. \\
     ${\mathbf \dim\Delta=4}.$ In the situation Case 1 the desired inequality is the following
     $$(F^{\Delta}+1)\cdot \ldots\cdot (F^{\Delta}+5)\leq 5!\cdot f_4^{\Delta}\cdot (F^{\Delta}+5-2^4+1+1)\cdot (F^{\Delta}-8)\cdot (F^{\Delta}-7)\cdot (F^{\Delta}-6).$$
     Inequality (\ref{ordering}) yields $F^{\Delta}+5\leq 16\cdot f_4^{\Delta}+15\leq 24\cdot f_4^{\Delta}$ if $f_4^{\Delta}\geq 2$, 
     i.e. $\Delta$ is not the $4$-simplex. It then suffices to show 
     $$(F^{\Delta}+1)\cdot(F^{\Delta}+2)\cdot (F^{\Delta}+3)\cdot (F^{\Delta}+4)\leq 5\cdot f_4^{\Delta}\cdot(F^{\Delta}-9)\cdot 
     (F^{\Delta}-8)\cdot (F^{\Delta}-7)\cdot (F^{\Delta}-6),$$
     which is equivalent to 
     $$4\cdot(F^{\Delta})^4-160\cdot (F^{\Delta})^3+1640\cdot (F^{\Delta})^2-8300\cdot F^{\Delta}+15096\geq 0.$$
     From $F^{\Delta}\geq 2^5-1+2^4-5=42$ for $f_4^{\Delta}\geq 2$ we deduce that $4\cdot F^{\Delta}\geq 160$ and $1640\cdot F^{\Delta}\geq 8300$. 
     This finally implies the claim.
     If $\Delta$ is the $4$-simplex then again the result follows from \cite{NS}.
     Now assume the situation of Case 2. If $\Delta$ is the boundary of the $5$-simplex then the assertion follows by Gorenstein-ness from 
     \cite{NS}. If $\Delta$ is not the boundary of the $5$-simplex then $f_4^\Delta \geq 7$. Here Kruskal-Katona theorem implies
      $f_3^\Delta , f_1^\Delta \geq 19, f_2^\Delta \geq 26, f0^\Delta \geq 7$. Thus $F^\Delta \geq 74$. 
      We we have to show that
      $$(F^{\Delta}+1)\cdot\ldots\cdot(F^{\Delta}+5)\leq 5\cdot f_4^{\Delta}\cdot(F^{\Delta}+5-2^5+2+1)\cdot (F^{\Delta}-23)\cdot \ldots\cdot (F^{\Delta}-21).$$
      By Inequality (\ref{facenumbers}) it holds that $F^{\Delta}+5\leq 16\cdot f_4^{\Delta}$. It thus suffices to show that 
      $$(F^{\Delta}+1)\cdot\ldots\cdot (F^{\Delta}+4)\leq \frac{15}{2} \cdot(F^{\Delta}-24)\cdot (F^{\Delta}-23)\cdot (F^{\Delta}-22)\cdot (F^{\Delta}-21).$$
%      i.e. 
%      $$(F^{\Delta})^4-1370\cdot (F^{\Delta})^3+45455\cdot (F^{\Delta})^2-681850\cdot F^{\Delta}+3825312 \geq 0.$$
      The latter inequality is true for $F^{\Delta}\geq 61$. Hence we are done. \\
      ${\mathbf \dim\Delta=5}.$ In the situation of Case 1 we have to show
      $$(F^{\Delta}+1)\cdot\ldots\cdot(F^{\Delta}+6)\leq 6!\cdot f_5^{\Delta}\cdot (F^{\Delta}+6-2^5+1+1)\cdot (F^{\Delta}-23)\cdot \ldots\cdot(F^{\Delta}-20).$$
      Inequality (\ref{ordering}) implies $F^{\Delta}+6\leq 32\cdot f_5^{\Delta}+31\leq 48\cdot f_5^{\Delta}$ if $f_5^{\Delta}\geq 2$, i.e. $\Delta$ is not the $5$-simplex. It then suffices to show that
      $$(F^{\Delta}+1)\cdot \ldots\cdot (F^{\Delta}+5)\leq 15\cdot (F^{\Delta}-24)\cdot \ldots\cdot (F^{\Delta}-20).$$
%      i.e.
%      $$14\cdot (F^{\Delta})^5-1665\cdot (F^{\Delta})^4+72440\cdot (F^{\Delta})^3-1592475\cdot (F^{\Delta})^2+17460086\cdot (F^{\Delta})-76507320\geq 0.$$
      The latter inequality is true for $F^{\Delta}\geq 57$. By $\dim\Delta=5$ it follows that $F^{\Delta}\geq 2^6-1-6=57$ which implies the assertion.\\
      In the situation of Case 2 we have to show that 
      $$(F^{\Delta}+1)\cdot \ldots\cdot(F^{\Delta}+6)\leq 6!\cdot f_5^{\Delta}\cdot(F^{\Delta}+6-2^6+2+1)\cdot (F^{\Delta}-54)\cdot \ldots\cdot (F^{\Delta}-51).$$
      By Inequality (\ref{facenumbers}) we know that $F^{\Delta}+6\leq 2^5\cdot f_5^{\Delta}$. Therefore, it suffices to show that
      $$(F^{\Delta}+1)\cdot \ldots\cdot (F^{\Delta}+5)\leq 22,5\cdot (F^{\Delta}-55)\cdot \ldots\cdot (F^{\Delta}-51).$$
%      i.e.
%      $$43\cdot (F^{\Delta})^5-11955\cdot (F^{\Delta})^4+1263655\cdot (F^{\Delta})^3-66959325\cdot (F^{\Delta})^2+1173461782\cdot F^{\Delta}-18785309640\geq 0.$$
      The above inequality is satisfied for $F^{\Delta}\geq 118$. Since by Inequality (\ref{face_number}) $F^{\Delta}\geq 2^7-1-1-6=120$ this implies 
      the claim. This concludes the proof of the lower bound part of the Multiplicity Conjecture if every $(\dim\Delta-1)$-dimensional face of $\Delta$ 
      lies in at least two facets of $\Delta$.\\

      ${\mathbf \dim\Delta=6}.$ We only need to consider the situation of Case 1. We have to show that
      $$(F^{\Delta}+1)\cdot \ldots\cdot (F^{\Delta}+7)\leq 7!\cdot f_6^{\Delta}\cdot (F^{\Delta}+7-2^6+1+1)\cdot(F^{\Delta}-54)\cdot \ldots\cdot (F^{\Delta}-50).$$
      From Inequality (\ref{ordering}) we deduce that $F^{\Delta}+7\leq 2^6-1+2^6\cdot f_6^{\Delta}=64\cdot f_6^{\Delta}+63\leq 96\cdot f_6^{\Delta}$ 
      if $f_6^{\Delta}\geq 2$, i.e. $\Delta$ is not the $6$-simplex. It then suffices to show that       
      $$(F^{\Delta}+1)\cdot \ldots\cdot (F^{\Delta}+6)\leq 52,5\cdot(F^{\Delta}-55)\cdot \ldots\cdot (F^{\Delta}-50).$$
%      i.e. $103\cdot (F^{\Delta})^6-33117\cdot (F^{\Delta})^5+4339825\cdot (F^{\Delta})^4-303685095\cdot (F^{\Delta})^3+11949944272\cdot (F^{\Delta})^2-250736330628\cdot F^{\Delta}+2191619428560\geq 0$.
      This finally concludes the proof of the lower bound of the Multiplicity Conjecture since the latter inequality is satisfied for 
       $F^{\Delta}\geq 113$ and since from $\dim\Delta=6$ we deduce that $F^{\Delta}\geq 2^7-1-7=120$.\\

      \noindent {\sf Cohen-Macaulay Case and Equality:} 
	   It remains to study the equality case when $k[\sd(\Delta)]$ is 
	   Cohen-Macaulay. By Section \ref{Cohen-Macaulay-ness} we
	   know that $k[\sd(\Delta)]$ is Cohen-Macaulay if and only
	   if $k[\Delta]$ is. 
	   Assume $\dim (\Delta) \geq 2$. 
	   From the fact that $\Delta$ is Cohen-Macaulay we infer that
	   either $\Delta$ is the $2$-simplex or $\Delta$ contains 
	   two $2$-dimensional faces that intersect along a 
	   $1$-dimensional face. If $\Delta$ is a $2$-simplex then
	   by inspection one sees that the resolution of
	   $k[\sd(\Delta)]$ is pure -- indeed linear -- and
	   satisfies the Multiplicity Conjecture with equality.
	   Hence after relabeling we may assume that $\Delta$
	   contains the face $\{1,2,3\}$ and the face $\{1,2,4\}$. 
	   We show that in this case the minimal free 
	   resolution of $k[\sd(\Delta)]$ is never pure and
	   the inequality always strict.
	   The restriction of $\sd(\Delta)$ to the vertices
	   $\{1\}$, $\{1,2,3\}$, $\{2\}$ and $\{1,2,4\}$ yields
	   a $4$-gon which shows by Corollary \ref{Betti_0} that 
           $\beta_{2,4} \neq 0$. On the
	   other hand the restriction of $\sd(\Delta)$ to the
	   vertices $\{1\}$, $\{2\}$, $\{3\}$ yields three
	   isolated points which then again by Corollary \ref{Betti_0}
           shows that $\beta_{2,3} \neq 0$.
	   Thus the minimal free resolution can never be pure. 
	   But our reasoning also implies that $M_2 \geq 4$. Since
	   our estimates only use $M_2 \geq 3$ this then shows that the
	   inequality is strict.
	   The case $\dim(\sd(\Delta)) = 0$ was already covered by the arguments
           preceding the whole proof.
	   Hence it remains to
	   consider $\dim (\Delta) = \dim(\sd(\Delta)) = 1$.
           By Theorem 4.3 from \cite{NS} it follows that equality
	   implies pureness of the minimal free resolution
	   for $k[\sd(\Delta)]$.
           Assume the minimal free resolution of $k[\sd(\Delta)]$
	   is pure. Let us first treat the case $\widetilde{H}_1(\Delta;k)= 0$.            
           Then Lemma \ref{lemma: Betti_2} 
	   along with the pureness of the minimal free resolution
	   imply $M_i=i+1$ for $1\leq i\leq f_0^{\Delta}+f_1^{\Delta}-2$.
	   Since $\Delta$ is a Cohen-Macaulay complex and 
           $\widetilde{H}_1(\Delta;k)=0$
	   the complex $\Delta$ is an acyclic connected graph and hence a tree. 
           Therefore, $f_0^{\Delta}=f_1^{\Delta}+1$.
	   Thus 
	   \begin{eqnarray*}		\frac{\prod_{i=1}^{f_0^{\Delta}+f_1^{\Delta}-2}M_i}{(f_0^{\Delta}+f_1^{\Delta}-2)!}&=&\frac{\prod_{i=1}^{f_0^{\Delta}+f_1^{\Delta}-2}(i+1)}{(f_0^{\Delta}+f_1^{\Delta}-2)!}\\
		&=&\frac{(f_0^{\Delta}+f_1^{\Delta}-1)!}{(f_0^{\Delta}+f_1^{\Delta}-2)!}=f_0^{\Delta}+f_1^{\Delta}-1=2f_1^{\Delta}.
				\end{eqnarray*}
	Since $e(k[\sd(\Delta)])=(\dim \Delta+1) \cdot f_1^{\Delta}=2f_1^{\Delta}$ the Multiplicity Conjecture is
	satisfied with equality.\\
	Let now $\widetilde{H}_1(\Delta;k)\neq 0$. We are going to show that if the minimal
	free resolution for $k[\sd(\Delta)]$ is pure, then $\Delta$ is an 
        $f_0^{\Delta}$-gon.
	In this case $k[\sd(\Delta)]$ satisfies the Multiplicity Conjecture with equality.
	By assumption $\Delta$ -- considered as a graph -- contains at least one cycle. Assume there are $\geq 2$ cycles.
        In the sequel we identify a cycle in $\Delta$ with the associated
        $1$-dimensional subcomplex of $\Delta$. 
	Let $\sigma_1, \sigma_2$ be two cycles of $\Delta$ such that $|\mathring{\sigma}_1|, |\mathring{\sigma}_2|$
	are minimal and $|\mathring{\sigma}_1|\leq |\mathring{\sigma}_2|$. If $|\mathring{\sigma}_1|< |\mathring{\sigma}_2|$,
	choose two vertices $v_1,v_2$ from $\sigma_2$ such that $\sd(\Delta)_{\mathring{\sigma}_2\setminus \{\{v_1\},\{v_2\}\}}$
	consists of two connected components $\sigma_2^1,\sigma_2^2$. 
	Let $W\subseteq \sigma_2\setminus\{\{v_1\},\{v_2\}\}$ such that 
	$\sigma_2^1\cap W\neq \emptyset$, $\sigma_2^2\cap W\neq \emptyset$ and
	$|W|=|\mathring{\sigma}_1|-1$. By construction $\widetilde{H}_0(\sd(\Delta)_W;k)\neq 0$. Hence by Corollary \ref{Betti_0} $\beta_{|W|-1,|W|}=\beta_{|\mathring{\sigma}_1|-2,|\mathring{\sigma}_1|-1}\neq 0$. Since
	 $\widetilde{H}_1(\sd(\Delta)_{\mathring{\sigma}_1};k)\neq 0$ it follows from Corollary \ref{Betti_0} that
	$\beta_{|\mathring{\sigma}_1|-2,|\mathring{\sigma}_1|}\neq 0$. Thus the minimal free resolution 
	of $k[\sd(\Delta)]$ cannot be pure in this case.\\
	If $|\mathring{\sigma}_1|=|\mathring{\sigma}_2|$, choose $v_1,v_2$ as in the above case and 
	in addition let $v_3$ be a vertex of $\sigma_1$ such that $v_1\neq v_3$ and $v_2\neq v_3$.
	Then, $\sd(\Delta)_{\mathring{\sigma}_2\setminus\{\{v_1\},\{v_2\}\}\cup\{\{v_3\}\}}$ has at least 
	two connected components which implies by $|\mathring{\sigma}_1|=|\mathring{\sigma}_2|$ and Corollary \ref{Betti_0} that
	$\beta_{|\mathring{\sigma}_1|-2,|\mathring{\sigma}_1|-1}\neq 0$. Thus, the minimal free resolution
	of $k[\sd(\Delta)]$ cannot be pure if $\Delta$ contains more than one cycle.\\
	Let $\sigma$ be the only cycle of $\Delta$. It remains to show that $\Delta$ cannot have 
	vertices respectively edges which do not lie in $\sigma$. If $v$ is a vertex of $\Delta$ such
	that $v\notin \sigma$ we can assume -- by $\Delta$ being connected --
	that there exists $w\in \sigma$ such that $\{v,w\}\in \Delta$. 
	$\sd(\Delta)_{\mathring{\sigma}\setminus\{\{v\},\{v,w\}\}\cup\{\{w\}\}}$ then consists of two 
	connected components what implies that $\beta_{|\mathring{\sigma}|-2,|\mathring{\sigma}|-1}\neq 0$.
	Since $\widetilde{H}_1(\sd(\Delta)_{\mathring{\sigma}};k)\neq 0$ the minimal free resolution
	of $k[\sd(\Delta)]$ cannot be pure. This completes the proof.
    \end{proof}
\section*{Acknowledgment}
  The authors are grateful to J\"urgen Herzog for suggesting the study of 
  the Multiplicity Conjecture for
  Stanley-Reisner rings of barycentric subdivisions and to Tim R\"omer for suggesting to study the equality situation.

  \end{document}